\documentclass[12pt]{article}
\usepackage{fullpage}
\usepackage{epsfig}
\usepackage{color} 
\usepackage{amscd}
\usepackage{latexsym}
\usepackage{amsfonts}
\usepackage{graphicx}
\usepackage{xypic}
\usepackage{amssymb, amsmath}
\numberwithin{equation}{section}
\usepackage[ps2pdf,                
bookmarks=true,
bookmarksnumbered=true,
hypertexnames=false,
breaklinks=true,
pdfborder={0 0 10}]{hyperref}

\def\Bbb{\mathbb}
\def\b{\mathbb}
\def\l{\label}
\def\bl{\bullet}

\def\08{\{0,\infty\}}
\def\01{\{0,1\}}

\def\e08{\t E_0\oplus \t_\infty}

\def\bx{\bar X}
\def\x8{X_\infty}
\def\os{\overset}
\def\0{\{0\}}

\def\t{\tilde}
\def\op{\oplus}
\def\os{\overset}

\def\d{\delta}

\def\im{\imath}
\def\jm{\jmath}

\def\08{\{0,\infty\}}
\def\e08{\t E_0\oplus \t_\infty}

\def\m1{(m-1)[2m-1]}

\def\ba{\mathbb A}
\def\bz{\mathbb Z}
\def\bzt{{\mathbb Z}_{tr}}
\def\bp{\mathbb P}
\def\cp{{\cal P}}

\def\r{\ref}
\def\er{\eqref}

\def\bean{\begin{eqnarray}}
\def\eean{\end{eqnarray}}
\def\bea{\begin{eqnarray*}}
\def\eea{\end{eqnarray*}}

\def\ben{\begin{enumerate}}
\def\een{\end{enumerate}}

\def\nm{\nonumber}

\def\lra{\longrightarrow}
\def\os{\overset}
\def\tvi{\widetilde{V_i}}
\def\io{\os\cong\to}

\newtheorem{thm}{Theorem}[section]
\newtheorem{defn}[thm]{Definition}
\newtheorem{main}[thm]{Main Theorem}
\newtheorem{prop}[thm]{Proposition}
\newtheorem{cor}[thm]{Corollary}
\newtheorem{lem}[thm]{Lemma}

\newtheorem{rmk}[thm]{Remark}

\newenvironment{pfpz}{\medskip \noindent
{\it Proof of Theorem \ref{pz}.}}{\hfill $\Box$
\\}
\newenvironment{pfmn}{\medskip \noindent
{\it Proof of Theorem \ref{mn}.}}{\hfill $\Box$
\\}

\newenvironment{proof}{\medskip \noindent
{\it Proof.}}{\hfill $\Box$
\\}

\title{Karoubi's Construction for Motivic\\ Cohomology Operations}

\author{Zhaohu Nie}

\date{March 2006}

\begin{document}

\maketitle

\begin{abstract}

We use an analogue of Karoubi's construction \cite k in the motivic situation to give some cohomology operations in motivic cohomology. We prove many properties of these operations, and we show that they coincide, up to some nonzero constants, with the reduced power operations in motivic cohomology originally constructed by Voevodsky \cite v. The relation of our construction to Voevodsky's is, roughly speaking, that of a fixed point set to its associated homotopy fixed point set. 


\end{abstract}

\tableofcontents

\section{Introduction}\l{1}

We fix a base field $k$ and a prime $l$ different from $char(k)$, the characteristic of $k$. 
Voevodsky \cite {v} constructed the reduced power operations 
$$P^i: H^{*,*}(-,{\Bbb Z}/l)\to H^{*+2i(l-1),*+i(l-1)}(-,{\Bbb Z}/l),\ i\geq 0$$ 
in motivic cohomology following the classical construction of Steenrod \cite s (see also \cite {k0}) for the reduced power operations in singular cohomology. Karoubi \cite k constructed these operations in singular cohomology using another construction. It is the purpose of this paper to adopt Karoubi's construction in the motivic situation. Our main construction and results are as follows. 

We work in the ${\b A}^1$-homotopy category $H^{\ba^1}_\bullet(k)$ of Morel and Voevodsky (see Section 2 and also \cite {mv,icm98,dv}). In this category, motivic cohomology is represented by motivic Eilenberg-MacLane spaces (see (\ref{rep})), and so motivic cohomology operations are represented by morphisms between the corresponding motivic Eilenberg-MacLane spaces in $H^{\ba^1}_\bullet(k)$. 

Let $n\geq 1$ be an integer. Choose our model of the motivic Eilenberg-MacLane space $K(\bz/l,2n,n)$ as $z_{equi}^{\bz/l}(\ba^n,0)$ (see \er{em}). The $l$-th cup product power operation 
\bean
P:K(\bz/l,2n,n)\to K(\bz/l,2nl,nl)\l{0p}
\eean
is represented by (abusing the notation $P$)
\bean
P: z_{equi}^{\bz/l}(\ba^n,0)\to z_{equi}^{\bz/l}(\ba^{nl},0),\l{0p1}
\eean
which is the $l$-th fiber product power of cycles (see \er{P}).

The cyclic group $C_l$ acts on $\ba^{nl}$ by cyclically permuting the $l$ copies of $\ba^n$, and so acts on $z^{\bz/l}_{equi}(\ba^{nl},0)$ by functoriality. By its very construction, the image of $P$ in \er{0p1} lies in the fixed point set $z^{\bz/l}_{equi}(\ba^{nl},0)^{C_l}$ of invariant cycles. 
Call the factorization 
\bean\l{0sp}
{\cal P}: z_{equi}^{\bz/l}(\ba^n,0)\to z_{equi}^{\bz/l}(\ba^{nl},0)^{C_l}
\eean
the total reduced power operation (see \er{op}). We have the following main theorem.

\begin{main}\l{main} Suppose that the field $k$ satisfies some assumptions (detailed at the beginning of Section \ref{3.2}).
\begin{enumerate}

\item In the $\ba^1$-homotopy category $H^{\ba^1}_\bullet(k)$, there is a canonical isomorphism  
$$z_{equi}^{\bz/l}(\ba^{nl},0)^{C_l}\cong \prod_{i=0}^{n(l-1)} K(\bz/l,2n+2i,n+i)
                         \times \prod_{i=0}^{n(l-1)-1} K(\bz/l,2n+2i+1,n+i).$$

\item Denote the components of 
$\cp$ in \er{0sp}, under the canonical isomorphism of part 1, by 
$$D^i: K(\bz/l,2n,n)\to K(\bz/l,2n+2i,n+i),\ 0\leq i\leq n(l-1)$$
and 
$$E^i: K(\bz/l,2n,n)\to K(\bz/l,2n+2i+1,n+i),\ 0\leq i\leq n(l-1)-1.$$
Then 
$$E^i=\beta\circ D^i,\ 0\leq i\leq n(l-1)-1,$$ 
where $\beta: K(\bz/l,*,*)\to K(\bz/l,*+1,*)$ represents the Bockstein homomorphism associated to the exact sequence of coefficients $0\to \bz/l\os {\cdot l} \to \bz/l^2\to \bz/l\to 0$.  

\item $D^0=id: K(\bz/l,2n,n)\to K(\bz/l,2n,n)$, and $D^{n(l-1)}=P:K(\bz/l,2n,n)\to K(\bz/l,2nl,nl)$, where $P$ is the $l$-th cup product power \er{0p}. 

\item For $i$ not a multiple of $(l-1)$, $D^i=0$ and thus $E^i=0$ by the conclusion of part 2. Up to nonzero constants in $\bz/l$, 
$$D^{j(l-1)}\equiv P^j,\ 0\leq j\leq n,$$
where the $P^j$ are the reduced power operations of Voevodsky \cite {v}. 
\end{enumerate}
\end{main}

The proof of this main theorem occupies most of the paper. In Section \r 2 we give some preliminaries. The proof of part 1 is done in Sections \r{3.1}, \r{3.2} and \r{3.3} (see Theorem \ref{phomtype}), part 2 in Section \r{3.4} (see Theorem \ref{pbock}), part 3 in Section \r{3.5} (see Theorems \ref{p0} and \ref{lastp}) and part 4 in Section \r{4.2} (see Theorem \ref{tcomp}). 

As will be seen in Section \r{4.2}, the relation of our construction to Voevodsky's is, roughly speaking, that of a fixed point set to its associated homotopy fixed point set.

\medskip\noindent{\it Acknowledgements} This paper is based on the author's Ph.D. thesis. He thanks his advisor, H. Blaine Lawson, Jr. for his guidance, discussions and support. He is also glad to thank many people for their useful discussions during the elaboration of this paper, including M. Karoubi, P. Lima-Filho and P. dos Santos. He thanks V. Voevodsky for his encouragement 
and useful discussions. Special thanks are due to Christian Haesemeyer, who answered a lot of the author's questions, gave a lot of invaluable suggestions, and read various versions of this paper carefully and critically.  
The author heartily thanks him. 


\section{Some background}\l{2}

Throughout this paper, we work in both Voevodsky's triangulated category of motives $DM$ \cite {dm} and the ${\b A}^1$-homotopy category $H^{{\b A}^1}_\bullet(k)$ constructed by Morel and Voevodsky \cite {{mv},{icm98},dv} over our base field $k$. In this section we give a quick review of the two categories and some of their properties relevant to this paper. The reader is referred to the above mentioned references for further details. For later use, we emphasize the concept of fundamental classes, and we prove a version of the projective bundle formula for motives with compact support. 

Let $Sm/k$ be the category of smooth schemes over $k$. There are several useful Grothendieck topologies on $Sm/k$, and the one most relevant to this paper is the Nisnevich topology (see \cite[Definition 3.1.2]{mv}). By definition, a presheaf of sets (abelian groups) on $Sm/k$ is a contravariant functor from $Sm/k$ to the category of sets (resp. abelian groups). A Nisnevich sheaf is a presheaf which satisfies the usual sheaf axioms for the Nisnevich topology. We denote the category of Nisnevich sheaves (of sets) on $Sm/k$ by $Shv_{Nis}(Sm/k)$. 

Let $SmCor(k)$ be the category of finite correspondences, whose objects are smooth schemes over $k$ and whose morphisms from $X$ to $Y$ are finite correspondences, i.e. algebraic cycles on $X\times Y$ which are finite and surjective over a component of $X$. $SmCor(k)$ is an additive category. An additive contravariant functor from $SmCor(k)$ to the category of abelian groups is called a presheaf with transfers. There is a natural functor $Sm/k\to SmCor(k)$, which is the identity on objects and which sends a morphism $f:X\to Y$ to its graph $\Gamma_f\subset X\times Y$. Therefore we can regard a presheaf with transfers as a presheaf of abelian groups on $Sm/k$. If the corresponding presheaf happens to be a Nisnevich sheaf, we will then call the presheaf with transfers a Nisnevich sheaf with transfers. We denote the category of Nisnevich sheaves with transfers by $Shv_{Nis}(SmCor(k))$. 

For a scheme $X$ of finite type over $k$, there are two naturally associated Nisnevich sheaves with transfers: ${\b Z}_{tr}(X)$ and $z_{equi}(X,0)$. By definition, for any smooth scheme $U\in Sm/k$, ${\b Z}_{tr}(X)(U)$ is the abelian group of algebraic cycles on $U\times X$ which are finite and surjective over a component of $U$, and $z_{equi}(X,0)(U)$ is the abelian group of algebraic cycles on $U\times X$ which are equidimensional of relative dimension 0 (quasi-finite) and dominant over a component of $U$. It can be checked that both ${\b Z}_{tr}(X)$ and $z_{equi}(X,0)$ are Nisnevich sheaves with transfers. When $X$ is compact, clearly $\bzt(X)=z_{equi}(X,0)$. Both ${\b Z}_{tr}(X)$ and $z_{equi}(X,0)$, being Nisnevich sheaves, have well defined pullbacks for an arbitrary morphism $f:U'\to U$, and we denote the pullbacks by $Cycl(f)$.  

Actually, one has the following two functors
\bean
{\b Z}_{tr}(-): Sch/k\to Shv_{Nis}(SmCor(k));\nm\\
z_{equi}(-,0): Sch^{prop}/k\to Shv_{Nis}(SmCor(k)),\l{zcov}
\eean
where $Sch/k$ is the category of schemes of finite type over $k$ and morphisms, and $Sch^{prop}/k$ is the category of schemes of finite type over $k$ and proper morphisms.

In general, for a scheme $X\in Sch/k$, one defines a Nisnevich sheaf with transfers $z_{equi}(X,r)$ for any integer $r\geq 0$ as follows. For a smooth scheme $U\in Sm/k$, 
$z_{equi}(X,r)(U)$ is defined to be the abelian group of algebraic cycles on $U\times X$ which are equidimensional of relative dimension $r$ and dominant over a component of $U$. If $f:Y\to X$ is a flat equidimensional morphism of relative dimension $n$, 
then flat pullback of cycles gives a morphism 
\bean\l{zpull}
f^*:z_{equi}(X,0)\to z_{equi}(Y,n).
\eean
In the particular case that $\im:U\to X$ is an open embedding, we usually denote the pullback $\im^*$ by 
\bean\l{zpull0}
r: z_{equi}(X,0)\to z_{equi}(U,0)
\eean
and call it the restriction map. 

In this paper we will mainly deal with $z_{equi}(X,0)$. Therefore for simplicity of notation we write $z(X)$ for it. For an abelian group $A$, $z_{equi}^A(X,0)=z_{equi}(X,0)\otimes A$ is the sheaf of equidimensional cycles with coefficients in $A$, and again we write $z^A(X)$ for short. 

$Shv_{Nis}(SmCor(k))$ is an abelian category \cite[Theorem 3.1.4]{dm}. We will denote its derived category of bounded from above complexes of Nisnevich sheaves with transfers by \\
$D^-(Shv_{Nis}(SmCor(k)))$.

A presheaf with transfers $F$ is called homotopy invariant if for any $X\in Sm/k$, the natural map $F(X)\to F(X\times \ba^1)$ induced by the projection $X\times \ba^1\to X$ is an isomorphism. A Nisnevich sheaf with transfers is called homotopy invariant if it is homotopy invariant as a presheaf with transfers. 

The category of (effective) motivic complexes $DM^{eff}_-(k)$, $DM$ for short, is defined to be the full subcategory of $D^-(Shv_{Nis}(SmCor(k)))$, which consists of all complexes with homotopy invariant cohomology sheaves. $DM$ is a tensor triangulated category. 

For a presheaf with transfers ${\cal F}$, one defines its singular complex $C_*{\cal F}$ as follows. Let $\Delta^\bullet$ be the standard cosimplicial objects in $Sm/k$ with $\Delta^n=Spec(k[x_0,x_1,\cdots,x_n]/\sum x_i=1)$ and the standard coface and codegeneracy maps. One defines $C_n{\cal F}$ by specifying 
$${ C}_n {\cal F}(U)={\cal F}(\Delta^n\times U)$$ 
for a smooth scheme $U$. The differential in $C_*{\cal F}$ is defined to be the alternating sum of the face maps. The cohomology presheaves of $C_*{\cal F}$ are homotopy invariant, and so are their Nisnevich sheafification \cite[Lemma 3.2.1]{dm}. Therefore $C_*{\cal F}\in DM$ for a Nisnevich sheaf with transfers ${\cal F}$. 

In particular, for a scheme $X$ of finite type over $k$, one defines its motive $M(X)$ and motive with compact support $M^c(X)$ in $DM$ as 
$$M(X)= C_* {\b Z}_{tr}(X),\ M^c(X)= C_* z(X).$$
When $X$ is compact, clearly $M(X)=M^c(X)$.
In view of \er{zcov}, one has the following functorial properties:
\bean
M(-): Sch/k\to DM;\ M^c(-): Sch^{prop}/k\to DM. \l{mcov}
\eean
For a proper morphism $f:Y\to X$, we denote the induced morphism on motives with compact support by 
$$f_*:M^c(Y)\to M^c(X).$$

There are several equivalent definitions of the motivic complexes $\bz(n)$ (whose shifts are also called Tate motives when considered as objects in $DM$). The one most relevant to us is the following natural isomorphism \cite{allchar}
\bean\l{an}
a_n:\bz(n)[2n]\os\cong\to C_*z(\ba^n)=M^c(\ba^n).
\eean

In $DM$, motivic cohomology is represented by Tate motives: for a scheme $X$ of finite type over $k$  
\bean
H^{j,i}(X)= Hom_{DM}(M(X), {\b Z}(i)[j]). \l{dmrep}
\eean
Dually, Borel-Moore motivic homology is defined by 
\bean
H^{BM}_{j,i}(X)=Hom_{DM}(\bz(i)[j],M^c(X)). \l{bmh}
\eean

In the remaining part of this section before we start to discuss the $\ba^1$-homotopy category, we suppose that $k$ is a field which admits resolution of singularities, since the results below are only proved at the moment under this condition. 

If $k$ is a field which admits resolution of singularities, then motives with compact support have the following tensor structure for product of schemes (see \cite[Proposition 4.1.7]{dm}):
\bean
M^c(X)\otimes M^c(Y)\os\cong\to M^c(X\times Y). \label{tensor}
\eean



The following localization distinguished triangle is of particular importance to our later computations. 
\begin{prop} \cite[Proposition 4.1.5]{dm} Suppose that $k$ is a field which admits resolution of singularities, $X$ is a scheme of finite type over $k$ and $Z$ is a closed subscheme of $X$. Then we have a canonical distinguished triangle of the form
\bean
M^c(Z)\to M^c(X)\to M^c(X-Z)\to M^c(Z)[1]. \label{ldt}
\eean
\end{prop}

The proof of this proposition uses a natural isomorphism induced from restriction
$$r:C_*(z(X)/z(Z))\os\cong\to C_*z(X-Z).$$

In particular, suppose that a non-compact variety $X$ has a compactification $\bar X$ with complement $X_\infty$:
$$X\subset \bar X,\ X_\infty=\bar X- X.$$
Then there is a natural isomorphism induced from restriction
\bean\l{locmap}
r_X: C_*(\bzt(\bx)/\bzt(\x8))\io C_*z(X)=M^c(X).
\l{rx}
\eean

The following result is a basic tool for us to produce maps to motives with compact support. 

\begin{prop}\l{direction}\cite[Corollary 4.2.4]{dm} Suppose that $k$ is a field which admits resolution of singularities. Let $f:Y\to X$ be a flat equidimensional morphism of relative dimension $n$ of schemes of finite type over $k$. Then there is a canonical morphism in $DM$ of the form: 
$$f^*: M^{c}(X)(n)[2n]\to M^c(Y),$$
and these morphisms satisfy all the standard properties of the contravariant functoriality of algebraic cycles. 
\end{prop}

\begin{proof} We apply the $C_*$ functor to \er{zpull} and get a morphism
$$M^c(X)=C_*z(X)\to C_*z_{equi}(Y,n).$$
\cite[Proposition 4.2.8]{dm} asserts a canonical isomorphism
\bean\l{zn?}
C_*z_{equi}(Y,n)= \underline{Hom}(\bz(n)[2n],M^c(Y)),
\eean
where the right hand side is an internal $Hom$ object in $DM$. The composition of the above two morphism is 
$$M^c(X)\to \underline{Hom}(\bz(n)[2n],M^c(Y)),$$
whose adjoint
$$f^*: M^{c}(X)(n)[2n]\to M^c(Y)$$
is the map that we are after. 
\end{proof}

Suppose that $X$ is a scheme of finite type over $k$ of dimension $n$. Consider the structure morphism $p:X\to Spec(k)$. Then using Proposition \ref{direction}, one defines the fundamental class of $X$ as the morphism 
\bean\l{fundclass}
cl_X:=p^*: \bz(n)[2n]\to M^c(X),
\eean
which in view of \er{bmh} represents a Borel-Moore homology class in $H^{BM}_{2n,n}(X)$. 
More generally, if $\jm:Y\to X$ is a subvariety of dimension $m$, then the following composition 
\bean\l{cly}
cl_Y:\bz(m)[2m]\to M^c(Y)\to M^c(X),
\eean
where the first arrow is the fundamental class of $Y$ and the second arrow is induced by the closed embedding $\jm$ using \er{mcov}, represents the fundamental class of $Y$ in $X$ as a Borel-Moore homology class in $H^{BM}_{2m,m}(X)$.  

As a special case of \er{fundclass}, one has the following fundamental class of $\ba^n$
\bean\l{cla}
cl_{\ba^n}:\bz(n)[2n]\os\cong \to M^c(\ba^n),
\eean
which is an isomorphism and which is the same as \er{an}. 

We have the following two simple results about fundamental classes. 

\begin{lem}\l{l-rest} Suppose that both $Y$ and $X$ are schemes of dimension $n$ over $k$, and one has a flat morphism $f:Y\to X$ of relative dimension 0  (e.g. an open embedding), i.e. one has the following diagram  
$$\xymatrix{
Y\ar[r]^f\ar[rd]_{p_Y} & X\ar[d]^{p_X}\\
 & k
,}$$
where the $p$'s are structure morphisms. 
Then there is the following commutative diagram 
$$\xymatrix{
M^c(Y) & M^c(X)\ar[l]_{f^*}\\
& \bz(n)[2n]\ar[u]_{p_X^*=cl_X}\ar[ul]^{p_Y^*=cl_Y}
.}$$
\end{lem}
\begin{proof} This follows from the contravariant functoriality of the map in Proposition \ref{direction} for flat morphisms. 
\end{proof}

\begin{lem}\l{l-push} Suppose that both $Y$ and $X$ are schemes of dimension $n$ over $k$, and one has a proper morphism $f:Y\to X$, i.e. one has the following diagram  
$$\xymatrix{
Y\ar[r]^f\ar[rd]_{p_Y} & X\ar[d]^{p_X}\\
 & k
,}$$
where the $p$'s are structure morphisms. 
Suppose that the degree of $f:Y\to X$ is $d$, which is the degree of function field extension. Then
there is the following commutative diagram 
$$\xymatrix{
M^c(Y)\ar[r]^{f_*} & M^c(X)\\
\bz(n)[2n]\ar[u]^{p_Y^*=cl_Y}\ar[r]^{\cdot d}& \bz(n)[2n]\ar[u]_{p_X^*=cl_X}
,}$$
where $\cdot d$ is multiplication by $d$. 
\end{lem}

\begin{proof} Under the proper pushforward of cycles by $f$, the cycle $Y$ is mapped to the $d$ multiple of the cycle $X$. 
\end{proof}

Now let's introduce a version of the projective bundle formula for motives with compact support. Suppose that $\pi:E\to X$ is a vector bundle of rank $n$. In view of Proposition \ref{direction}, one has a natural morphism 
\bean\l{pi}
\pi^*: M^c(X)(n)[2n]\to M^c(E).
\eean

The map $\pi^*$ in \er{pi} has the following simple functorial properties. 

\begin{lem}\l{pullbdl} Let $f:Y\to X$ be a morphism of schemes, and $\pi:E\to X$ a vector bundle over $X$ of rank $n$. Let $f^*(E)$ be the pullback vector bundle on $Y$, i.e. we have the following Cartesian diagram
\bean
\xymatrix{
f^*(E)\ar[r]^g\ar[d]^\pi & E\ar[d]^\pi\\
Y\ar[r]^f & X.
}\l{f*E}
\eean

If $f$ is proper, and therefore so is $g$, then one has the following commutative diagram
\bean
\xymatrix{
M^c(f^*(E))\ar[r]^{g_*} & M^c(E)\\
M^c(Y)(n)[2n]\ar[r]^{f_*\otimes id}\ar[u]^{\pi^*} & M^c(X)(n)[2n]\ar[u]^{\pi^*},
}\l{pb_*}
\eean
where the vertical morphisms are as in \er{pi}, $f_*$ and $g_*$ are the induced maps on motives with compact support by the properness of $f$ and $g$. 

On the other hand, if $f$ is flat of relative dimension 0, and therefore so is $g$, then one has the following commutative diagram
\bean
\xymatrix{
M^c(f^*(E)) & M^c(E)\ar[l]_{g^*}\\
M^c(Y)(n)[2n]\ar[u]^{\pi^*} & M^c(X)(n)[2n]\ar[l]_{f^*\otimes id}\ar[u]^{\pi^*},
}\l{pb^*}
\eean
where $f^*$ and $g^*$ are the induced maps on motives with compact support by the flatness of $f$ and $g$. 

\end{lem}

\begin{proof} One has the following commutative diagram on the cycle level
$$\xymatrix{
z_{equi}(f^*(E),n)\ar[r]^{g_*} & z_{equi}(E,n)\\
z(Y)\ar[r]^{f_*}\ar[u]^{\pi^*} & z(X)\ar[u]^{\pi^*} ,
}$$
since \er{f*E} is Cartesian. Apply the $C_*$ functor, and we conclude \er{pb_*} by the naturality of \er{zn?} for proper pushforward (see the proof of Proposition \ref{direction}). 


\er{pb^*} holds by the contravariant functoriality of the map in Proposition \ref{direction} for flat morphisms. 
\end{proof}


We have the following version of the projective bundle formula. 

\begin{thm}\l{proj-mc} Suppose that $k$ admits resolution of singularities. The map $\pi^*$ in \er{pi} is an isomorphism. 
\end{thm}

\begin{proof} Let $\im:U\to X$ be an open subscheme such that $E|_U$ is trivial. Let $\jm:Y=X-U\to X$ be the closed complement. Consider the following diagram of localization distinguished triangles \er{ldt}:
$$\xymatrix{
M^c(E|_Y)\ar[r] & M^c(E)\ar[r] & M^c(E|_U) \\
M^c(Y)(n)[2n]\ar[r]^{\jm\otimes id}\ar[u]^{\pi^*} & M^c(X)(n)[2n]\ar[r]^{r\otimes id}\ar[u]^{\pi^*} & M^c(U)(n)[2n]\ar[u]^{\pi^*}.
}$$
The commutativity of the above diagram follows from Lemma \ref{pullbdl}. 

Since $E|_U\cong U\times \ba^n$ is trivial, the last vertical arrow 
$$\pi^*: M^c(U)(n)[2n]\to M^c(E|_U)$$
is an isomorphism by \er{cla} and \er{tensor}. Then an induction on the dimension of $X$ and the five lemma complete the proof. 
\end{proof}




Now let's start to review the ${\b A}^1$-homotopy theory $H^{{\b A}^1}_\bullet(k)$ (hence removing the assumption that $k$ admits resolution of singularities). We start with the category $sShv_{Nis,\bullet}(Sm/k)$ 
of simplicial Nisnevich sheaves of pointed sets over $Sm/k$. A Nisnevich sheaf is always identified to its associated constant simplicial sheaf, i.e. all the terms of the simplicial sheaf are the given sheaf and all the face and degeneracy maps are the identities. 
$H^{{\b A}^1}_\bullet (k)$ is the successive localization of $sShv_{Nis,\bullet}(Sm/k)$ with respect to the classes of simplicial weak equivalences and ${\b A}^1$-weak equivalences (see \cite{mv,icm98,dv} for more details). 

For an abelian group $A$, 
let $K(A,j,i)$ be the simplicial abelian sheaf corresponding to the complex $A(i)[j]=\bz(i)[j]\otimes A$ under the Dold-Kan correspondence functor $K$. Considered as simplicial sheaf of pointed sets, it defines an object of $H^{{\b A}^1}_\bullet (k)$. It is a  motivic Eilenberg-MacLane space, since it represents the corresponding motivic cohomology \cite[Theorem 2.3.1]{dv} as 
\bean
Hom_{H^{{\b A}^1}_\bullet (k)}(X_+,K(A,j,i))=H^{j,i}(X,A), \label{rep}
\eean
where $X$ is a smooth scheme and $X_+$ is regarded in $H^{{\b A}^1}_\bullet(k)$ as the constant simplicial representable sheaf associated to $X$ with a disjoint base point. 

We can choose a model of the motivic Eilenberg-MacLane space $K(A,2n,n)$ as 
$$K(A(n)[2n])=K(C_* z^A(\ba^n))$$
in view of \er{an}. In addition, by \cite[Lemma 2.5.2]{dv}, we have a natural $\ba^1$-weak equivalence 
$$K(z^A(\ba^{n}))\os\cong\to K(C_* z^A(\ba^{n})),$$ 
where $z^A(\ba^{n})$ on the left hand side is considered to be a complex concentrated at dimension 0. It is easy to see by the construction of $K$ in the Dold-Kan correspondence that $K(z^A(\ba^{n}))$ is the constant simplicial sheaf at $z^A(\ba^{n})$. Therefore we have the following nice model for $K(A,2n,n)$ as
\bean\l{em}
K(A,2n,n)=z^{A}({\b A}^n).
\eean


Because of the representability (\ref{rep}) of motivic cohomology by Eilenberg-MacLane spaces in the ${\b A}^1$-homotopy category $H^{{\b A}^1}_\bullet (k)$, motivic cohomology operations are the same as morphisms between corresponding Eilenberg-MacLane spaces in this category. 
 
Actually our construction will map the above model $z^{\bz/l}({\b A}^n)$ of $K({\Bbb Z}/l,2n,n)$ to various Eilenberg-MacLane spaces of other dimensions by performing some geometric constructions on it. 

\section{The total reduced power operation}\l{3.1}

In this section, we work with an arbitrary coefficient ring $R$, and the ones of interest to us are $\bz/l$ and $\bz$. We omit the coefficient ring from our notation for simplicity. 

Let's start with our model of the motivic Eilenberg-MacLane space $K(2n,n)$ as $z(\ba^n)$ \er{em}. Consider the $l$-th cup product power operation
$$P:K(2n,n)\to K(2nl,nl),$$
which is represented by the $l$-th fiber product power of cycles
\bean
P:z(\ba^n)\to z(\ba^{nl}).\label{P}
\eean
In more detail, for a smooth scheme $U$ and a cycle $Z\in z(\ba^n)(U)$, write $Z=\sum n_i Z_i$, where $n_i\in R$ and the $Z_i$ are closed and irreducible subschemes of $U\times \ba^n$, which are equidimensional of relative dimension 0 over $U$. Then 
$$P(Z)=\sum n_{i_1}\cdots n_{i_l} Z_{i_1}\times_U\cdots \times_U Z_{i_l},$$
where $Z_{i_1}\times_U\cdots \times_U Z_{i_l}$ is the fiber product of $Z_{i_1},\cdots,Z_{i_l}$ over $U$. It is easy to see that $P(Z)\in z(\ba^{nl})(U)$.  Note that the map $P$ is not linear, i.e. $P$ is a morphism in the category of sheaves of sets but not of $R$ modules. Figure \ref{fig-P} illustrates the map $P$ \er{P} when $l=2$ and $Z$ is an irreducible subscheme.

\begin{figure}
\centering
\input{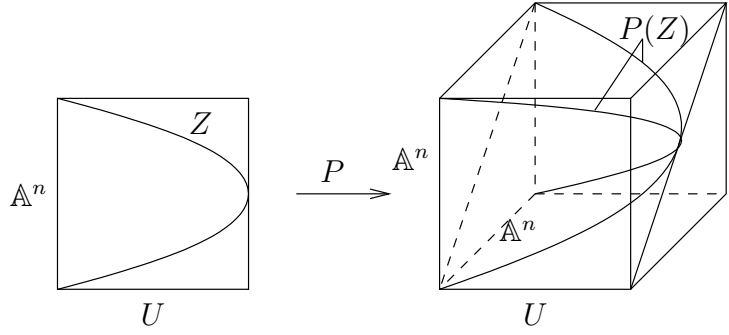}
\caption{The power operation $P$}
\l{fig-P}
\end{figure}


The symmetric group $S_l$ acts on $\ba^{nl}$ by permuting the $l$ copies of $\ba^n$, and in particular the cyclic group $C_l$ acts on $\ba^{nl}$ by cyclically permuting the $l$ copies. This naturally induces an action of $S_l$, and in particular of $C_l$, on $z(\ba^{nl})$ by functoriality. The important thing for us is that the above map (\ref{P}) factorizes as
\bean
P:z(\ba^n)\overset {{\cal P}}{\lra} z(\ba^{nl})^{C_l}\os{j}{\lra} z(\ba^{nl}),\l{factorization}
\eean
where $(-)^{C_l}$ stands for the fixed point set under the $C_l$ action on the sheaf level, and $j$ is the natural inclusion map. More precisely, $z(\ba^{nl})^{C_l}$ is the Nisnevich sheaf whose value on a smooth scheme $U$ is the subgroup of cycles in $z(\ba^{nl})(U)$ which are invariant under the $C_l$ action on $U\times \ba^{nl}$ as the product of the trivial action on $U$ and the cyclic permutation on $\ba^{nl}$. Concretely an algebraic cycle $Z\in z(\ba^{nl})(U)$ is invariant if and only if all its fibers over points of $U$, as linear combinations of points in $\ba^{nl}$, are invariant under the $C_l$ action on $\ba^{nl}$. From this we see that $z(\ba^{nl})^{C_l}$ is still a presheaf (pullbacks of invariant cycles are invariant). Actually it is a Nisnevich sheaf with transfers as one can easily check. 
Furthermore the image of $\cp$ lies in the similarly defined fixed point set $z(\ba^{nl})^{S_l}$ by the whole symmetric group. 


We call the map 
\bean
{\cal P}:z(\ba^n)\to z(\ba^{nl})^{C_l} \l{op}
\eean
the total reduced power operation. Our task is to analyze the homotopy type of $z(\ba^{nl})^{C_l}$ in the $\ba^1$-homotopy category $H^{{\b A}^1}_\bullet (k)$.

One important feature of the $C_l$ action on $\ba^{nl}$ when $l$ is a prime is its semi-freeness. The diagonal of $\ba^{nl}$, denoted by $\ba^n$, is the fixed point set, and on the complement, $\ba^{nl}-\ba^n$, the action is free.


Using the functorialities \er{zcov} and \er{zpull0}, one has the following exact sequence of sheaves
\bean
0\lra z(\ba^n)\os \Delta \lra z(\ba^{nl})^{C_l}\os r \lra z(\ba^{nl}-\ba^n)^{C_l},\label{es}
\eean
where the inclusion $\Delta$ is induced by the closed embedding of the diagonal $\ba^n$ in $\ba^{nl}$, and the restriction $r$ is induced by the open embedding of the complement $\ba^{nl}-\ba^n$ in $\ba^{nl}$, as illustrated in Figure \ref{fig-es}.

\begin{figure}
\centering
\input{r.eps}
\caption{The exact sequence \er{es}}
\l{fig-es}
\end{figure}

\begin{prop}\l{pfdt}
Let $k$ be a field which admits resolution of singularities. We have a distinguished triangle in $DM$
\bean
C_*z(\ba^n)\os\Delta\to C_*z(\ba^{nl})^{C_l}\os r \to C_*z(\ba^{nl}-\ba^n)^{C_l} \os \delta \to C_*z(\ba^n)[1].\label{fdt}
\eean
\end{prop}

\begin{proof} The proof is an analogue of \cite [Theorem 5.11] {fv}, except that we have to take care of the group action here. We write out all the detail for the convenience of the reader. The exact sequence \er{es} clearly induces the following diagram with the row as a distinguished triangle in $DM$
$$\xymatrix{
C_*z(\ba^n)\ar[r]^{C_*\Delta}& C_*z(\ba^{nl})^{C_l}\ar[r]& C_*(z(\ba^{nl})^{C_l}/z(\ba^n))\ar[r]\ar[d]^{C_*r} & C_*z(\ba^n)[1].\\
 & &C_* z(\ba^{nl}-\ba^n)^{C_l}& 
}$$
To prove the proposition, we only need to prove that $C_*r$ is a quasi-isomorphism, or equivalently $C_*(z(\ba^{nl}-\ba^n)^{C_l}/z(\ba^{nl})^{C_l})$ is acyclic, i.e. quasi-isomorphic to zero. First let's recall the following criterion for acyclicity.

\begin{prop}\cite [Theorem 4.1.2] {dm} Let $k$ be a field which admits resolution of singularities, and $F$ be a presheaf with transfers on $Sm/k$ such that for any smooth scheme $U$ over $k$ and a section $\phi\in F(U)$ there is a proper birational morphism $p:U'\to U$ with $U'$ smooth and 
$F(p)(\phi)=0$. Then the complex $C_*(F)$ is acyclic. 
\end{prop}

In view of this criterion, it suffices to show that for any smooth scheme $U$ over $k$ and a closed integral subscheme $Z$ in $z(\ba^{nl}-\ba^n)^{C_l}(U)$ that there is a blow-up $p:U'\to U$ with $U'$ smooth 
such that $Cycl(p)(Z)$, the pullback cycle of $Z$ by $p$, lies in the image of $z(\ba^{nl})^{C_l}(U')$. 


By the platification theorem (see \cite[Theorem 2.2.3]{sv2} and \cite{flat}), there is a blow-up $p:U'\to U$ such that the closure of $Cycl(p)(Z)\in z(\ba^{nl}-\ba^n)(U')$ lies in the image of $z(\ba^{nl})(U')$. 
It is clear that the closure $\overline{Cycl(p)(Z)}\in z(\ba^{nl})^{C_l}(U')$ since $Cycl(p)(Z)\in z(\ba^{nl}-\ba^n)^{C_l}(U')$ and the fibers added when taking the closure are all supported on the diagonal $\ba^n$, which is the fixed point set under the $C_l$ action. 
\end{proof}

Due to the freeness of the action of $C_l$ on $\ba^{nl}-\ba^n$, we have the following canonical mutually inverse isomorphisms 
\bean
z(\ba^{nl}-\ba^n)^{C_l}\os {\widetilde{\pi_*}}\lra z((\ba^{nl}-\ba^n)/C_l)\os {\widetilde{\pi^*}} \lra z(\ba^{nl}-\ba^n)^{C_l}, \label{free}
\eean
where $\widetilde{\pi^*}$ is the flat pullback whose image is obviously invariant, and $\widetilde{\pi_*}$ is the reduced pushforward, i.e. for $Z\in z(\ba^{nl}-\ba^n)^{C_l}(U)$, $\widetilde{\pi_*}(Z)$ is the unique cycle in $z((\ba^{nl}-\ba^n)/C_l)(U)$ whose flat pullback is $Z$. The situation when $l=2$ is illustrated in Figure \ref{fig-free}.

\begin{figure}
\centering
\input{pi.eps}
\caption{The isomorphism \er{free}}
\l{fig-free}
\end{figure}

Following the language of Karoubi \cite k, we call $(\ba^{nl}-\ba^n)/C_l$ the $l$-th normalized cyclic product of $\ba^n$, and sometimes denote it by $CP^+_l(\ba^n)$. 

In the next section, we will compute the motive type of 
$$M^c((\ba^{nl}-\ba^n)/C_l)=C_*z((\ba^{nl}-\ba^n)/C_l).$$ 
Together with (\ref{free}) and (\ref{fdt}), we are able to analyze the motive type of $C_*z(\ba^{nl})^{C_l}$. The motive type of $C_*z^{\bz/l}(\ba^{nl})^{C_l}$ for $\bz/l$ coefficients can then be easily deduced. This then determines the homotopy type of $z^{\bz/l}(\ba^{nl})^{C_l}$. We will see that $z^{\bz/l}(\ba^{nl})^{C_l}$ has the desired homotopy type to be the target of the total reduced power operation \er{op}. 

\section{Lens spaces}\l{3.2}

From now on we assume that $k$ admits resolution of singularities to be able to use the tools described in Section 2 and Proposition \ref{pfdt}. 
This condition at the moment means that $char(k)=0$. 
We furthermore assume that $k$ has a primitive $l$-th root of unity $\zeta$ to simplify our exposition. Note that under the condition that $l\neq char(k)$, this assumption doesn't constitute a real restriction, and can always be achieved by passing to a separable extension of $k$ of degree prime to $l$. Transfer arguments show that this doesn't affect considerations for motivic cohomology with $\bz/l$ coefficients, which is the main subject in this paper. 
Under this assumption we have an isomorphism between the cyclic group $C_l$ and $\mu_l$, the group of $l$-th roots of unity in $k$, which sends $1\in \bz/l=C_l$ to $\zeta\in \mu_l$. Note that $\mu_l$ has a standard representation on ${\b A}^1$, denoted by $\rho$, where $\zeta\in \mu_l$ acts as $z\to \zeta\cdot z$ for $z\in \ba^1$.


Under the assumption of the existence of $\zeta$, it is clear that if one chooses a suitable basis of ${\b A}^{nl}$, the action of $C_l$ on ${\b A}^{nl}$ by cyclically permuting the $l$ copies of ${\b A}^n$ can be written as a direct sum $id\oplus \rho\oplus \rho^2\oplus \cdots\oplus\rho^{l-1}$. Here $id$ is the trivial action, the action $\rho$ is interpreted through the isomorphism of $C_l$ to $\mu_l$, $\rho^i$ is the $i$-th tensor product of $\rho$ where $\zeta\in \mu_l$ acts as $z\to \zeta^i\cdot z$, and for notational simplicity we have written $\rho^i$ for $\os {n} {\overbrace{\rho^i\oplus\cdots\oplus\rho^i}}$, which acts on a copy of ${\b A}^n$. 



Therefore under this new basis, we can rewrite our normalized cyclic product as
$$CP^+_l({\b A}^n)=({\b A}^{nl}-{\b A}^n)/C_l={\b A}^n\times ({\b A}^{n(l-1)}-\{0\})/C_l,$$ 
where $C_l$ acts freely on ${\b A}^{n(l-1)}-\{0\}$.


By our model of the Tate motives \er{an} and the tensor structure of motives with compact support (\ref{tensor}), we have the following sequence of isomorphisms 
\begin{eqnarray}
&& M^c(({\Bbb A}^{n(l-1)}-\{0\})/C_l)\otimes {\Bbb Z}(n)[2n]\nm\\
&\os\cong\to& M^c(({\Bbb A}^{n(l-1)}-\{0\})/C_l)\otimes M^c({\Bbb A}^{n})\nm\\
&\io& M^c({\b A}^n\times ({\b A}^{n(l-1)}-\{0\})/C_l)=M^c(({\Bbb A}^{nl}-{\Bbb A}^{n})/C_l).\label{tensor1}
\end{eqnarray}

Now we want to decide the motive type of $M^c(({\Bbb A}^{n(l-1)}-\{0\})/C_l)$. This is a special case of the following general considerations. Suppose that $\mu_l$ acts on ${\Bbb A}^m$ as a direct sum of nontrivial irreducible representations $\rho^{a_1}\oplus\cdots\oplus\rho^{a_m}$ with $1\leq a_i\leq l-1$ for $1\leq i\leq m$, where $m\geq 1$ is a general dimension. Then $(\ba^m-\{0\})/\mu_l$ is a lens space of dimension $m$, and we denote it by $L_m$. (The actual weights $(a_1,\cdots,a_m)$ are not important in our later computations.) We want to determine the motive type of $M^c(L_m)$. 

The space $\ba^m$ has the following natural filtration
$$\ba^1\subset\cdots\subset \ba^m,$$
where $\ba^i$ is the subspace of points with the last $(m-i)$ coordinates zero. This induces a filtration on $\ba^m-\{0\}$
$$\ba^1-\{0\}\subset\cdots\subset \ba^m-\{0\},$$
and a filtration on $L_m$
$$L_1\subset\cdots\subset L_m,$$
where $L_i=(\ba^i-\{0\})/\mu_l$ is the skeletal lens space of dimension $i$. Let 
$$\pi: \ba^i-\{0\}\to L_i$$ 
be the natural projection. 

We also denote the complement of $L_{m-1}$ in $L_m$ by $U_m$:
$$U_m:=L_m-L_{m-1}=(\ba^m-\ba^{m-1})/\mu_l
=(\ba^{m-1}\times (\ba^1-\{0\}))/\mu_l,$$
where $\ba^1$ is the last component of $\ba^m$. Since the action of $\mu_l$ on $\ba^1-\{0\}$ is free, $U_m$ is a vector bundle of rank $(m-1)$ over $(\ba^1-\{0\})/\mu_l$.
We will often use the notation $\ba^*:=\ba^1-\{0\}$. 

The main result of this section is the following theorem. 

\begin{thm}\l{pz} There is a canonical isomorphism
\bean\l{gm}
f_m:\bz(0)[1]\op\bigoplus_{i=1}^{m-1} \bz/l(i)[2i]\op \bz(m)[2m]\os\cong \to M^c(L_m).
\eean
\end{thm}

We will prove this theorem by induction. To study the case $m=1$, we first establish the following result. 

\begin{prop}\l{pa1*} There is a canonical isomorphism 
$$s_1\op cl: \bz(0)[1]\op \bz(1)[2]\os \cong \to M^c(\ba^1-\{0\}).$$
\end{prop}

\begin{proof} We define $cl: \bz(1)[2]\to M^c(\ba^1-\{0\})$ to be the fundamental class of $\ba^1-\{0\}$ \er{fundclass}. We now construct $s_1: \bz(0)[1]\to M^c(\ba^1-\{0\})$. 
One has the following two natural isomorphisms:
\ben
\item
\bea
\sigma_1:\bz(0)[1]
=C_*(\bzt(\ba^1)/\bzt(\{0,1\})).
\eea
Actually $\ba^1/\{0,1\}$ is one model of the simplicial circle $S^1_s$ in the homotopy category $H_\bullet^{\ba^1}(k)$. (Also the above isomorphism can be seen from the following distinguished triangle
$$
        C_*\bzt(\{0,1\})\to C_*\bzt(\ba^1)\to C_*(\bzt(\ba^1)/\bzt(\{0,1\}))\to C_*\bzt(\{0,1\})[1],
$$
which has the form 
$$\bz\op\bz\to \bz\to C_*(\bzt(\ba^1)/\bzt(\{0,1\}))\to \bz[1]\op\bz[1],$$
where the first map is $(a,b)\mapsto a+b$.)

\item Considering the natural compactification $\bp^1$ of $\ba^1-\{0\}$ with complement $\{0,\infty\}$ and using \er{rx}, one has a natural isomorphism 
$$r_{\ba^1-\{0\}}: 
C_*(\bzt(\bp^1)/\bzt(\{0,\infty\}))\io C_*z(\ba^1-\{0\})=M^c(\ba^1-\{0\}).$$

\een


Let's fix an automorphism
\bean\l{auto}
\phi:(\bp^1,\{0,1\})\to (\bp^1,\{0,\infty\});\ z\mapsto \frac z {z-1}.
\eean
(In homogeneous coordinates, the map $\phi$ is $[x_0,x_1]\mapsto [x_1-x_0, x_1]$.)

Let $\imath_1': (\ba^1,\{0,1\})\to (\bp^1,\{0,1\})$ be the natural inclusion, and
$$\imath_1=\phi\circ \im_1': (\ba^1,\{0,1\})\to (\bp^1,\{0,1\})\to (\bp^1,\{0,\infty\})$$
the composition. 
Abusing notation, we also denote the induced morphism after applying $C_*\bzt(-)$ by $\im_1$:
$$\im_1:C_*(\bzt(\ba^1)/\bzt(\{0,1\}))\to C_*(\bzt(\bp^1)/\bzt(\{0,\infty\}).$$ 
Then we define $s_1$ as the composition 
\bean
s_1=r_{\ba^1-\{0\}}\circ \im_1\circ \sigma_1&:& \bz(0)[1]=C_*(\bzt(\ba^1)/\bzt(\{0,1\}))\nm\\
&\to& C_*(\bzt(\bp^1)/\bzt(\{0,\infty\})\to M^c(\ba^1-\{0\}).\l{o1}
\eean

Now we want to prove the morphism 
$$s_1\op cl: \bz(0)[1]\op \bz(1)[2]\to M^c(\ba^1-\{0\})$$
constructed above is an isomorphism. 

Consider the following diagram of distinguished triangles
$$\xymatrix{
\bz(0)[0]\ar[r]^0\ar[d]^{cl_{\{0\}}} & \bz(1)[2]\ar[r]\ar[d]^{cl_{\ba^1}} & \bz(0)[1]\op\bz(1)[2]\ar[r]\ar[d]^{s_1\op cl} &  \bz(0)[1]\ar[d]^{cl_{\{0\}}[1]}\\
M^c(\{0\})\ar[r] &  M^c(\ba^1)\ar[r] & M^c(\ba^1-\{0\})\ar[r]^\delta & M^c(\{0\})[1].
}
$$
Here the second row is the localization distinguished triangle \eqref{ldt} associated to the embeddings 
$$\{0\}\os{\rm closed}\lra \ba^1\os{\rm open}\longleftarrow \ba^1-\{0\}.$$
The maps in the first row are the natural ones: 0, inclusion and projection. 
The two first vertical maps are induced by the fundamental classes of the point $0\in \ba^1$ and $\ba^1$, and they are known to be isomorphisms \er{cla}.

Now we want to show that the diagram commutes. The commutativity of the first square is clear, and the commutativity of the second one follows from Lemma \ref{l-rest}. Thus we only need to prove the commutativity of the third square on the $\bz(0)[1]$ factor. 

By our construction, 
\bea
\delta\circ s_1&:& \bz(0)[1]\to C_*(\bzt(\bp^1)/\bzt(\{0,\infty\})\\
&\to& {\rm ker}(C_*\bzt(\{0,\infty\})[1] \to C_*\bzt(\bp^1)[1]).
\eea
Since any two points on $\bp^1$ are rationally equivalent, the kernel in the last map is equivalent to $C_*\bzt(\{0\})[1]=M^c(\{0\})[1]$, and the natural map from $\bz(0)[1]$ to it is given by the fundamental class of $\{0\}$ shifted by [1]. 

Now we conclude by the five lemma. 
\end{proof}

\begin{rmk} More generally, the automorphism in \er{auto} has the form 
\bean\l{ak*}
z\mapsto \frac {az} {z-1},
\eean
where $a\in k^*$. In view of Proposition \ref{pa1*}, the Hom set $Hom(\bz(0)[1],M^c(\ba^1-\{0\}))$ has a factor ${Hom}(\bz(0)[1],\bz(1)[2])=H^{1,1}(k)=k^*$. The $a$ in \er{ak*} realizes this factor. 
\end{rmk}




We now prove the $m=1$ case of Theorem \ref{pz} in the following lemma. 

\begin{lem}\l{m=1} There is a natural isomorphism 
$$f_1={}_0f_1\op {}_1f_1: \bz(0)[1]\op \bz(1)[2]\os\cong \to M^c(L_1).$$
Furthermore, one has the following two commutative diagrams
\bean
\xymatrix{
\bz(0)[1]\op \bz(1)[2]\ar[r]^{id\op(\cdot l)}\ar[d]^{s_1\op cl} & \bz(0)[1]\op \bz(1)[2]\ar[d]^{{}_0f_1\op {}_1f_1}\\
M^c(\ba^1-\{0\})\ar[r]^{\pi_*} & M^c(L_1);
}\l{onel}
\eean
\bean
\xymatrix{
\bz(0)[1]\op \bz(1)[2]\ar[d]^{s_1\op cl} & \bz(0)[1]\op \bz(1)[2]\ar[l]_{(\cdot l)\op id}\ar[d]^{{}_0f_1\op {}_1f_1}\\
M^c(\ba^1-\{0\}) & M^c(L_1)\ar[l]_{\pi^*}.
}\l{lone}
\eean
\end{lem}

\begin{proof} We define 
\bean
{}_0f_1:=\pi_*\circ s_1: \bz(0)[1]\to M^c(\ba^1-\{0\})\to M^c(L_1)
\l{0g1}
\eean
and
$${}_1f_1:=cl:\bz(1)[2]\to M^c(L_1)$$
as the fundamental class of $L_1$. 

$L_1=(\ba^1-\{0\})/\mu_l$ can be canonically identified to $\ba^1-\{0\}$. Under this identification, the above defined ${}_0f_1$ and ${}_1f_1$ become $s_1$ and $cl$ in Proposition \ref{pa1*}. Therefore $f_1$ is an isomorphism. 

The commutativity of the diagram \er{onel} on the $\bz(0)[1]$ factor is by the definition of ${}_0f_1$ \er{0g1}, and that on the $\bz(1)[2]$ factor is by Lemma \ref{l-push} since $\pi:\ba^1-\{0\}\to L_1$ is proper of degree $l$. 

The commutativity of the diagram \er{lone} follows from \er{onel} since one has $\pi_*\circ \pi^*=(\cdot l)\op (\cdot l).$


\end{proof}

Since $U_m$ is a vector bundle of rank $(m-1)$ over $(\ba^1-\{0\})/\mu_l$, one has the following canonical isomorphism 
\bean
s\op cl:&& {\b Z}(m-1)[2m-1]\oplus {\b Z}(m)[2m]= (\bz(0)[1]\op\bz(1)[2])(m-1)[2m-2]\nm\\
\os\cong\to&& M^c((\ba^1-\{0\})/\mu_l)(m-1)[2m-2]\os\cong\to M^c(U_m),
\l{erum}
\eean
where $\bz(0)[1]\op\bz(1)[2]\to M^c((\ba^1-\{0\})/\mu_l)$ is an analogue of $f_1$ in Lemma \ref{m=1}, and the last isomorphism is by the projective bundle formula, Theorem \ref{proj-mc}.

Now let's start to prove Theorem \ref{pz}. 

\begin{pfpz} For notational simplicity, we often denote the left hand side of \er{gm} by $F_m$. We denote the $i$-th component of $f_m$ on the $i$-th summand of $F_m$ (in the order in \er{gm} with the first one as the 0-th) by ${}_if_m$. 

For a general dimension $m$, we always define the last component ${}_mf_m$ by the fundamental class of $L_m$ \er{fundclass}, i.e.
\bean
{}_mf_m=cl_{L_m}:\bz(m)[2m]\to M^c(L_m).
\l{mgm}
\eean

To define the other components of $f_m$, 
one actually needs an induction procedure, which we start now. 

The case when $m=1$ has been discussed in Lemma \ref{m=1}. 

Now we assume that
\bean\l{indhyp}
f_{m-1}:\bz(0)[1]\op\bigoplus_{i=1}^{m-2} \bz/l(i)[2i]\op \bz(m-1)[2m-2]\to M^c(L_{m-1})
\eean
has been canonically defined and is an isomorphism. 

Using the natural morphism $\jmath:M^c(L_{m-1})\to M^c(L_m)$ induced by the closed embedding of $L_{m-1}$ in $L_m$, we define 
$${}_0f_m:=\jmath\circ {}_0f_{m-1}: \bz(0)[1]\to M^c(L_{m-1})\to M^c(L_m),$$
\bean\l{lowerones}
{}_if_m:=\jmath\circ {}_if_{m-1}: \bz/l(i)[2i]\to M^c(L_{m-1})\to M^c(L_m),\ 1\leq i\leq m-2.
\eean
With ${}_mf_m$ defined by the fundamental class \er{mgm}, we only need to define 
\bean\l{m-1fm}
{}_{m-1}f_m:\bz/l(m-1)[2m-2]\to M^c(L_m).
\eean

To this end, we need to use the following two distinguished triangles:
$$\bz(m-1)[2m-2]\os{\cdot l}\to \bz(m-1)[2m-2]\os\tau\to \bz/l(m-1)[2m-2]\os\beta\to \bz(m-1)[2m-1],$$
which is the Bockstein distinguished triangle, and 
$$M^c(L_{m-1})\os \jmath\to M^c(L_m)\os r\to M^c(U_m)\os \d \to M^c(L_{m-1})[1],$$
which is the localization distinguished triangle \er{ldt} associated to 
$$\xymatrix{
L_{m-1}\ar[r]^{\rm closed} & L_m & U_m\ar[l]_{\rm open}
.}$$ 

We consider the following diagram of long exact sequences of $Hom$ groups in the category $DM$
(to fit the rather big diagram in the paper, we use the following shorthand notation: $H$ for $Hom$, $\{\}$ for $(m-1)[2m-2]$, and  ${}^cL_{m-1}$ for $M^c(L_{m-1})$ etc. (apology for the notation)):
\bean\xymatrix{
\bullet\ar[r]^<<<<<<<{\cdot l}
& H(\bz\{\}[1],{}^cL_{m-1})\ar[r]^{\beta^*}\ar[d]^{\jmath_*} & H(\bz/l\{\},{}^cL_{m-1})\ar[r]^{\tau^*}\ar[d]^{\jmath_*} & H(\bz\{\},{}^cL_{m-1})  \ar[d]^{\jmath_*}\ar[r]^>>>>>>>{\cdot l} & 
\bl\\
\bl\ar[r]^<<<<<<<{\cdot l}
& H(\bz\{\}[1],{}^cL_m)\ar[r]^{\beta^*}\ar[d]^{r_*} & H(\bz/l\{\},{}^cL_m)\ar[r]^{\tau^*}\ar[d]^{r_*} & H(\bz\{\},{}^cL_m)\ar[d]^{r_*}\ar[r]^>>>>>>>{\cdot l} & 
\bl\\
\bl\ar[r]^<<<<<<<{\cdot l}
& H(\bz\{\}[1],{}^cU_m)\ar[r]^{\beta^*}\ar[d]^{\d_*} & H(\bz/l\{\},{}^cU_m)\ar[r]^{\tau^*}\ar[d]^{\d_*} & H(\bz\{\},{}^cU_m)\ar[d]^{\d_*}\ar[r]^>>>>>>>{\cdot l} & 
\bl\\
\bl\ar[r]^<<<<<<<{\cdot l} & H(\bz\{\}[1],{}^cL_{m-1}[1])\ar[r]^{\beta^*} & H(\bz/l\{\},{}^cL_{m-1}[1])\ar[r]^{\tau^*} & H(\bz\{\},{}^cL_{m-1}[1])\ar[r]^>>>>>>>{\cdot l} & 
\bl
}\l{chasing1}\eean

We have the following two explicit morphisms:
$$cl_{L_{m-1}}\in H(\bz\{\},{}^cL_{m-1})=Hom(\bz(m-1)[2m-2],M^c(L_{m-1})),$$
which is the fundamental class of $L_{m-1}$, and 
$$s\in H(\bz\{\}[1],{}^cU_m)=Hom(\bz(m-1)[2m-1],M^c(U_m)),$$
which is defined in \er{erum}. 

The following proposition gives the map ${}_{m-1}f_m$ in \er{m-1fm}.

\begin{prop}\l{m1gm} There is a unique element 
$${}_{m-1}f_m\in Hom(\bz/l(m-1)[2m-2],M^c(L_m)),$$
such that 
\bean
r_*({}_{m-1}f_m)=\beta^*(s).
\l{go}
\eean

In addition, one has 
\bean
\tau^*({}_{m-1}f_m)=\jmath_*(cl_{L_{m-1}}).
\l{gl}
\eean
\end{prop}

\begin{proof} First we want to show that 
\bean
\d_*(s)=l\cdot cl_{L_{m-1}}[1],
\l{dotl}
\eean
where $cl_{L_{m-1}}[1]\in H(\bz\{\}[1],{}^cL_{m-1})[1])$ 
is the fundamental class of $L_{m-1}$ (shifted by 1). This is equivalent to the commutativity of the following diagram

\bean\xymatrix{
\bz(m-1)[2m-1]\ar[r]^{\cdot l}\ar[d]^{s} & \bz(m-1)[2m-1]\ar[d]^{cl_{L_{m-1}}[1]}\\
M^c(U_m)\ar[r]^\d & M^c(L_{m-1})[1]
.}\l{delta1}\eean

To prove this, we first give the following general case of Proposition \ref{pa1*} for a general dimension $m$. 

\begin{prop}\l{am*} There is a natural isomorphism
$$s_m\op cl: \bz(0)[1]\op \bz(m)[2m]\os \cong \to M^c(\ba^m-\{0\}).$$
\end{prop}

\begin{proof} Define $cl:\bz(m)[2m]\to M^c(\ba^m-\{0\})$ to be the fundamental class of $\ba^m-\{0\}$ \er{fundclass}. We define $s_m$ as the following composition 
$$s_m:=\jm_m\circ s_1: \bz(0)[1]\to M^c(\ba^1-\{0\})\to M^c(\ba^m-\{0\}),$$
where $\jm_m$ is induced by the closed embedding $\ba^1-\{0\}\to \ba^m-\{0\}$ with the last $(m-1)$ coordinates zero. 

The proof of $s_m\op cl$ being an isomorphism is almost the same as that for Proposition \ref{pa1*}. We don't repeat here. 
\end{proof}

Consider the following diagram of embeddings
\bean
\xymatrix{
{\b A}^{m-1}-\{0\}\ar[r]^{\rm closed}\ar[d]^\pi &{\b A}^{m}-\{0\}\ar[d]^\pi &{\b A}^{m-1}\times {\b A}^*\ar[l]_{\rm open}\ar[d]^\pi \\
L_{m-1}\ar[r]^{{\rm closed}}& L_m & U_m\ar[l]_{\rm open},
}\label{sd}
\eean
where the vertical arrows are the quotient maps. 

By (\ref{ldt}), the diagram (\ref{sd}) gives a diagram of localization distinguished triangles
\bean
\xymatrix{
M^c({\b A}^{m-1}-\{0\})\ar[r]\ar[d]^{\pi_*} & M^c({\b A}^{m}-\{0\})\ar[r]\ar[d]^{\pi_*} & M^c({\b A}^{m-1}\times {\b A}^*)\ar[d]^{\pi_*} \ar[r]^{\d_1} & M^c({\b A}^{m-1}-\{0\})[1]\ar[d]^{\pi_*} \\
M^c(L_{m-1})\ar[r]& M^c(L_m)\ar[r]& M^c(U_m)\ar[r]^\d & M^c(L_{m-1})[1].
}\label{both}
\eean

We can analyze the first row in \er{both} by considering the following diagram of distinguished triangles
$$\xymatrix{
{\b Z}(0)[1]\oplus {\b Z}(m-1)[2m-2]\ar[r]\ar[d]^{s_{m-1}\op cl} & {\b Z}(0)[1]\oplus {\b Z}(m)[2m]\ar[r]\ar[d]^{s_m\op cl} & \\
M^c({\b A}^{m-1}-\{0\})\ar[r] & M^c({\b A}^{m}-\{0\})\ar[r] &  
}$$
\bean\xymatrix{
 \ar[r] & {\b Z}(m-1)[2m-1]\oplus {\b Z}(m)[2m] \ar[r] \ar[d]^{s\op cl} & {\b Z}(0)[2]\oplus {\b Z}(m-1)[2m-1]\ar[d]^{(s_{m-1}\op cl)[1]}\\
 \ar[r] & M^c({\b A}^{m-1}\times {\b A}^*)\ar[r]^{\d_1} & M^c(\ba^{m-1}-\{0\})[1]
.}\l{beforequotients}
\eean
Here the first row of distinguished triangle is a natural one consisting of inclusions and projections, and the vertical maps are those from Proposition \ref{am*} (the third one is by the tensor product structure \er{tensor} and Proposition \ref{pa1*}). It can be seen that the diagram \er{beforequotients} is commutative and is a natural isomorphism of distinguished triangles. 

Now consider the following diagram of the boundary maps of \er{both} with the natural maps from the Tate motive $\bz(m-1)[2m-2]$:
\bean
\xymatrix{
\bz(m-1)[2m-1]\ar[rr]^{id}\ar[rd]^{id}\ar[dd]^s & & \bz(m-1)[2m-1]\ar[rd]^{\cdot l}\ar'[d]^{cl[1]}[dd] & \\
 & \bz(m-1)[2m-1]\ar[rr]^(0.3){\cdot l}\ar[dd]^(0.3)s & & \bz(m-1)[2m-1]\ar[dd]^{cl[1]}\\
M^c(\ba^{m-1}\times \ba^*)\ar'[r][rr]^(0.3){\delta_1}\ar[rd]^{\pi_*}  & & M^c(\ba^{m-1}-\{0\})[1]\ar[rd]^{\pi_*} & \\
 & M^c(U_m)\ar[rr]^{\delta} & &  M^c(L_{m-1})[1]
}\l{cube1}
\eean


The left square of \er{cube1} is commutative by the following commutative diagram 
$$\xymatrix{
(\bz(0)[1])(m-1)[2m-2]\ar[r]^{id\otimes id}\ar[d]^{s_1\otimes id} & (\bz(0)[1])(m-1)[2m-2]\ar[d]^{{}_0f_1\otimes id}\\
M^c(\ba^*)(m-1)[2m-2]\ar[r]^{\pi_*\otimes id}\ar[d]^{p^*} & M^c(\ba^*/\mu_l)(m-1)[2m-2]\ar[d]^{p^*}\\
M^c(\ba^{m-1}\times \ba^*)\ar[r]^{\pi_*} &  M^c(U_m),
}$$
where the first square commutes by \er{onel}
, and the second by \er{pb_*}, 
since 
$$\xymatrix{
\ba^{m-1}\times \ba^*\ar[r]^\pi\ar[d]^p & U_m\ar[d]^p\\
\ba^*\ar[r]^\pi & \ba^*/\mu_l
}$$
is Cartesian. 

The right square in \er{cube1} is commutative by Lemma \ref{l-push}, and the back one is commutative by \er{beforequotients}. The top and bottom squares in \er{cube1} are commutative obviously. Therefore the front square, which is just diagram \er{delta1}, is commutative. (Here actually one considers the composition of the front square with the skew morphism $\bz(m-1)[2m-2]\os{id}\to\bz(m-1)[2m-2]$ in \er{cube1}.) Thus we have proved \er{dotl}. 

Therefore in diagram \er{chasing1}, 
$$\d_*(\beta^*s)=\beta^* (\d_* s)=\beta^* (l\cdot cl_{L_{m-1}}[1])=0.$$ 
Now by the induction hypothesis \er{indhyp},
one has the top term in the middle column of \er{chasing1}
$H(\bz/l\{\},{}^cL_{m-1})=Hom(\bz/l(m-1)[2m-2],M^c(L_{m-1}))=0.$

By the exactness of the second column in \er{chasing1}, there is a unique 
$${}_{m-1}f_m\in Hom(\bz/l(m-1)[2m-2],M^c(L_m)),$$ 
which satisfies \er{go}.

Now we want to prove \er{gl} for the above chosen ${}_{m-1}f_m$.
This is equivalent to the commutativity of the second square in the following diagram:
$$\xymatrix{ 
\bz(m-1)[2m-2]\ar[r]^{\cdot l}\ar[d]^{s[-1]} & \bz(m-1)[2m-2]\ar[r]^\tau\ar[d]^{cl_{L_{m-1}}} & \bz/l(m-1)[2m-2]\ar[r]^\beta\ar[d]^{{}_{m-1}f_m} & \bz(m-1)[2m-1]\ar[d]^{s}\\
M^c(U_m)[-1]\ar[r]^\d & M^c(L_{m-1}) \ar[r]^\jm & M^c(L_m)\ar[r]^r & M^c(U_m).
}$$

We know that the third square commutes by \er{go}. By the theory of triangulated categories, there exists a map $\psi:\bz(m-1)[2m-2]\to M^c(L_{m-1})$ for the second vertical arrow which makes both the first and the second squares commute.  If we use $cl_{L_{m-1}}$ for the second vertical arrow, the first square commutes by \er{dotl}. By our induction hypothesis \er{indhyp}, we see that 
$$Hom(\bz(m-1)[2m-2],M^c(L_{m-1}))=\bz$$ 
generated by $cl_{L_{m-1}}$. Therefore $\psi$ must be $cl_{L_{m-1}}$ to make the first square commute. Therefore $cl_{L_{m-1}}$ makes the second square commute, which proves \er{gl}. 
\end{proof}

\begin{rmk}
One can directly check that $l(\jm_*(cl_{L_{m-1}}))=0$, which is implied by \er{gl}. 
That is, we want to show the composition in the following diagram 
$$\xymatrix{
\bz(m-1)[2m-2]\ar[r]^(.6){cl_{L_{m-1}}} & M^c(L_{m-1})\ar[r]^\jm & M^c(L_m)\\
\bz(m-1)[2m-2]\ar[u]^{\cdot l} & &  
}$$
is zero. 

The above diagram can be embedded in the following one
$$\xymatrix{
\bz(m-1)[2m-2]\ar[r]^(.6){cl_{L_{m-1}}} & M^c(L_{m-1})\ar[r]^\jm & M^c(L_m)\\
\bz(m-1)[2m-2]\ar[r]^{cl
}\ar[u]^{\cdot l} & M^c(\ba^{m-1}-\{0\})\ar[u]^{\pi_*}\ar[r]^\jm & M^c(\ba^m-\{0\})\ar[u]^{\pi_*}
,}$$
where $cl$ is the fundamental class of $\ba^{m-1}-\{0\}$, and the $\jm$'s are induced by closed embeddings. The first square commutes by Lemma \ref{l-push}. 
The composition of the lower-right path is zero since it factorizes through $M^c(\ba^m-\{0\})$, to which any morphism from $\bz(m-1)[2m-2]$ is zero by Lemma \ref{am*} and dimensional consideration. 

\end{rmk}

Now we want to prove that the above defined morphism 
$$f_m:F_m\to M^c(L_m)$$
is an isomorphism. 

Consider the following diagram of distinguished triangles
\bean
\xymatrix{
F_{m-1}\ar[r]\ar[d]^{f_{m-1}} & F_m\ar[r]\ar[d]^{f_m} & \bz(m-1)[2m-1]\op \bz(m)[2m]\ar[r]\ar[d]^{s\op cl} & F_{m-1}[1]\ar[d]^{f_{m-1}[1]}\\
M^c(L_{m-1})\ar[r]^\jm & M^c(L_m)\ar[r]^r & M^c(U_m)\ar[r]^\d & M^c(L_{m-1})[1]
,}\l{gmiso}
\eean
where 
the third vertical morphism is that in \er{erum}. Besides natural inclusions and projections, the only nontrivial maps in the first row consist of 
$$\bz(m-1)[2m-2]\os{\tau}\to \bz/l(m-1)[2m-2]\os{\beta}\to \bz(m-1)[2m-1]\os{\cdot l}\to \bz(m-1)[2m-1].$$
Now we want to see that the diagram \er{gmiso} commutes. 

The commutativity of the first square on $\bz(m-1)[2m-2]$ is by \er{gl}, and on the other components by definition \er{lowerones}. 

The commutativity of the second square on $\bz(m)[2m]$ is by Lemma \ref{l-rest}, and on \\
$\bz/l(m-1)[2m-2]$ by \er{go}.

The commutativity of the third square on $\bz(m-1)[2m-1]$ has been shown in  \er{dotl}. 

Since $f_{m-1}$ and $(s\op cl)$ are isomorphisms by the induction hypothesis \er{indhyp} and \er{erum}, we conclude that $f_m$ is an isomorphism by the five lemma. 
\end{pfpz}

Theorem \ref{pz} and (\ref{tensor1}) give the following corollary concerning the motive with compact support of the normalized cyclic product $(\ba^{nl}-\ba^n)/C_l$.

\begin{cor} \l{pncp} In the triangulated category of motives DM, we have a natural isomorphism
\bean
h:\bz(n)[2n+1]\oplus\bigoplus_{i=n+1}^{nl-1}\bz/l(i)[2i]\oplus\bz(nl)[2nl]\io M^c(({\b A}^{nl}-{\b A}^n)/C_l).\l{ncp}
\eean
\end{cor}
 
We denote the first component of the above map by 
\bean\l{hn}
h_n:\bz(n)[2n+1]\to M^c(({\b A}^{nl}-{\b A}^n)/C_l).
\eean

We denote the left hand side of \er{ncp} by $H$. Let $h^{-1}: M^c(({\b A}^{nl}-{\b A}^n)/C_l)\to H$ be the inverse isomorphism. We define
\bean
h^{-1}_{>n}
:  M^c(({\b A}^{nl}-{\b A}^n)/C_l)\to H\to \bigoplus_{i=n+1}^{nl-1}\bz/l(i)[2i]\oplus\bz(nl)[2nl];\l{h-11}\\
h^{-1}_n: M^c(({\b A}^{nl}-{\b A}^n)/C_l)\to H\to \bz(n)[2n+1].\l{u1}
\eean
to be the compositions of $h^{-1}$ with the natural projections. 

\section{The total fixed point set}\l{3.3}

In this section, we want to study the motive type of $C_*z(\ba^{nl})^{C_l}$ using the distinguished triangle \er{fdt} and Corollary \ref{pncp}. The method we use is very similar to that in the proof of Theorem \ref{pz}. 

First of all, let's define a morphism 
\bean
g=(g_n,g_{>n}):C_*z(\ba^{nl})^{C_l}\to \bz/l(n)[2n]\op (\bigoplus_{i=n+1}^{nl-1}\bz/l(i)[2i]\op \bz(nl)[2nl]),
\l{definef}
\eean
where 
\bean
g_{>n}:=h_{>n}^{-1}\circ r:&& C_*z(\ba^{nl})^{C_l}\to C_*z(\ba^{nl}-\ba^n)^{C_l}\cong M^c((\ba^{nl}-\ba^n)/C_l)\nm\\
\to&& \bigoplus_{i=n+1}^{nl-1}\bz/l(i)[2i]\oplus\bz(nl)[2nl],\l{g>n}
\eean
and 
\bean\l{defgn}
g_n:C_*z(\ba^{nl})^{C_l}\to \bz/l(n)[2n]
\eean
is defined as follows. 

We consider the Bockstein distinguished triangle
$$\bz(n)[2n]\os{\cdot l}\to \bz(n)[2n]\os\tau\to \bz/l(n)[2n]\os\beta\to \bz(n)[2n+1]$$
and the distinguished triangle \er{fdt} 
$$C_*z(\ba^n)\os\Delta\to C_*z(\ba^{nl})^{C_l}\os r \to C_*z(\ba^{nl}-\ba^n)^{C_l}\os\delta\to C_*z(\ba^n)[1].$$
Again, we have the following diagram of long exact sequences of $Hom$ groups in $DM$, where we use shorthand notation such as $H$ for $Hom$, $\bz(n]$ for $\bz(n)[2n]$, $\ba^n$ for $C_*z(\ba^n)$, $\ba^{nlC}$ for $C_*z(\ba^{nl})^{C_l}$ and $\ba^{(nl-n)C}$ for $C_*z(\ba^{nl}-\ba^n)^{C_l}$. (Apology again for the notation.)

\bean\xymatrix{
\bullet\ar[r]^<<<<<<<{\cdot l}
 & H(\ba^{(nl-n)C},\bz(n])\ar[r]^{\tau_*}\ar[d]^{r^*} & H(\ba^{(nl-n)C},\bz/l(n])\ar[r]^{\beta_*}\ar[d]^{r^*} & H(\ba^{(nl-n)C},\bz(n][1])\ar[d]^{r^*}\ar[r]^>>>>>>>{\cdot l} & 
\bl\\
\bl\ar[r]^<<<<<<<{\cdot l}
 & H(\ba^{nlC},\bz(n])\ar[r]^{\tau_*}\ar[d]^{\Delta^*} & H(\ba^{nlC},\bz/l(n])\ar[r]^{\beta_*}\ar[d]^{\Delta^*} & H(\ba^{nlC},\bz(n][1])\ar[d]^{\Delta^*}\ar[r]^>>>>>>>{\cdot l} & 
\bl\\
\bl\ar[r]^<<<<<<<{\cdot l}
& H(\ba^n,\bz(n])\ar[r]^{\tau_*}\ar[d]^{\d^*} & H(\ba^n,\bz/l(n])\ar[r]^{\beta_*}\ar[d]^{\d^*} & H(\ba^n,\bz(n][1])\ar[d]^{\d^*}\ar[r]^>>>>>>>{\cdot l} & 
\bl\\
\bl\ar[r]^<<<<<<<{\cdot l} & H(\ba^{(nl-n)C}[-1],\bz(n])\ar[r]^{\tau_*} & H(\ba^{(nl-n)C}[-1],\bz/l(n])\ar[r]^{\beta_*} & H(\ba^{(nl-n)C}[-1],\bz(n][1])\ar[r]^>>>>>>>{\cdot l} & 
\bl
}\l{chasing2}\eean

We have the following two explicit morphisms:
$$cl_{\ba^n}^{-1}\in H(\ba^n,\bz(n])=Hom(C_*z(\ba^n),\bz(n)[2n]),$$
which is the inverse of \er{cla}, and 
$$h_n^{-1}\in H(\ba^{(nl-n)C},\bz(n][1])=Hom(C_*z(\ba^{nl}-\ba^n)^{C_l},\bz(n)[2n+1]),$$
which is as in \er{u1}. 

The following proposition gives the map $g_n$ in \er{defgn}, and it is very similar to Proposition \ref{m1gm}. 

\begin{prop}\l{ng} There is a unique element $g_n\in Hom(C_*z(\ba^{nl})^{C_l},\bz/l(n)[2n])$, such that 
\bean
\Delta^*(g_n)=\tau_*(cl_{\ba^n}^{-1}).
\l{fa}
\eean

In addition, one has
\bean
\beta_*(g_n)=r^*(h_n^{-1}).
\l{fu}
\eean
\end{prop}

\begin{proof} First we want to show that 
\bean
\d^*(cl_{\ba^n}^{-1})=l\cdot h_n^{-1}[-1], 
\l{ua}
\eean
where $h_n^{-1}[-1]\in H(\ba^{(nl-n)C}[-1],\bz(n])$ is what we had in \er{u1} (shifted by [-1]). 
This is equivalent to the commutativity of 
$$\xymatrix{
\bz(n)[2n]\ar[r]^(.6){\cdot l} & \bz(n)[2n]\\
C_*z(\ba^{nl}-\ba^n)^{C_l}[-1]\ar[r]^(.6)\d\ar[u]_{h_n^{-1}[-1]} & C_*z(\ba^n)\ar[u]_{cl_{\ba^n}^{-1}}
.}$$
Since $h_n:\bz(n)[2n+1]\to C_*z(\ba^{nl}-\ba^n)^{C_l}$ \er{hn} is a splitting of $h_n^{-1}$, the above commutativity, after shifting by $[1]$, is equivalent to that of the following diagram

\bean
\xymatrix{
\bz(n)[2n+1]\ar[r]^{\cdot l}\ar[d]^{h_n} & \bz(n)[2n+1]\ar[d]^{cl_{\ba^n}[1]}\\
C_*z(\ba^{nl}-\ba^n)^{C_l}\ar[r]^\d & C_*z(\ba^n)[1]
.}\l{delta2}
\eean

First consider the following diagram of distinguished triangles
\bean\xymatrix{
{\b Z}(n)[2n]\ar[r]^0\ar[d]^{cl} & {\b Z}(nl)[2nl]\ar[r]\ar[d]^{cl} & {\b Z}(n)[2n+1]\oplus {\b Z}(nl)[2nl] \ar[r]\ar[d]^{s\op cl} & {\b Z}(n)[2n+1]\ar[d]^{cl[1]}\\
C_*z(\ba^n)\ar[r] & C_*z(\ba^{nl})\ar[r] & C_*z(\ba^{nl}-\ba^n)\ar[r]^{\d_0} & C_*z(\ba^n)[1]
,}\l{whole}\eean
where the morphisms in the first row consist of 0, inclusion and projection, the second row is a natural localization distinguished triangle (\ref{ldt}), and the vertical maps are fundamental classes (the third one uses \er{tensor} and Proposition \ref{am*}). It can be seen that the above diagram is commutative and is a natural isomorphism of distinguished triangles. 

Now consider the following diagram of distinguished triangles
\bean
\xymatrix{
C_*z(\ba^n)\ar[r]\ar[d]^{id} & C_*z(\ba^{nl})^{C_l}\ar[r]\ar[d]^{C_*j} & C_*z(\ba^{nl}-\ba^n)^{C_l}\ar[r]^{\delta}\ar[d]^{C_*j} & C_*z(\ba^n)[1]\ar[d]^{id}\\
C_*z(\ba^n)\ar[r] & C_*z(\ba^{nl})\ar[r] & C_*z(\ba^{nl}-\ba^n)\ar[r]^{\delta_0} & C_*z(\ba^n)[1],
}\l{both1}
\eean
where the vertical arrows are induced by the inclusions $j$ of the fixed point set sheaves into the whole sheaves. Such a diagram exists due to the proof of Proposition \ref{pfdt}, which is analogous to the proof of the localization distinguished triangle \ref{ldt} (see \cite [Theorem 5.11] {fv}). Now let's try to understand the relevant parts of the third square in (\ref{both1}) concerning the factor $\bz(n)[2n+1]$.

We get the following diagram 
\bean\xymatrix{
\bz(n)[2n+1]\ar[rr]^{id}\ar[dd]^s & & \bz(n)[2n+1]\ar'[d]^{cl_{\ba^n}[1]}[dd] & \\
 & \bz(n)[2n+1]\ar[lu]^{\cdot l}\ar[rr]^(.3){\cdot l}\ar[dd]^(.3){h_n} & & \bz(n)[2n+1]\ar[dd]^{cl_{\ba^n}[1]}\ar[lu]^{id}\\
C_*z(\ba^{nl}-\ba^n)\ar'[r][rr]^(.3){\delta_0} & & C_*z(\ba^n)[1]& \\
 & C_*z(\ba^{nl}-\ba^n)^{C_l}\ar[rr]^{\delta}\ar[lu]^{C_*j} & &  C_*z(\ba^n)[1]\ar[lu]^{id}
}
\l{cube2}
\eean

\begin{lem}\l{leftsq} The left square in \er{cube2} is commutative. 
\end{lem}

\begin{proof} Under the identification \er{free}:
$$C_*z(\ba^{nl}-\ba^n)^{C_l}\cong C_*z((\ba^{nl}-\ba^n)/C_l)=M^c((\ba^{nl}-\ba^n)/C_l),$$
the map $C_*j$ in the left square of \er{cube2} is the same as 
$$\pi^*:M^c((\ba^{nl}-\ba^n)/C_l)\to M^c(\ba^{nl}-\ba^n).$$

The left square commutes by the following diagram 
$$\xymatrix{
(\bz(0)[1])(n)[2n]\ar[r]^{(\cdot l)\otimes id}\ar[d]^{({}_0f_1)\otimes id} & (\bz(0)[1])(n)[2n]\ar[d]^{s_{1}\otimes id}\\
M^c((\ba^1-\{0\})/\mu_l)(n)[2n]\ar[d]^{\jm\otimes cl_{\ba^n}}\ar[r]^{\pi^*} & M^c(\ba^1-\{0\})(n)[2n]\ar[d]^{\jm\otimes cl_{\ba^n}}\\
M^c((\ba^{n(l-1)}-\{0\})/C_l)\otimes M^c(\ba^n)\ar[d]^\cong\ar[r]^{\pi^*\otimes id} &  M^c(\ba^{n(l-1)}-\{0\})\otimes M^c(\ba^n)\ar[d]^\cong\\
M^c((\ba^{nl}-\ba^n)/C_l)\ar[r]^{\pi^*} & M^c(\ba^{nl}-\ba^n),
}$$
where the composition of the left column is $h_n$ \er{hn}, the composition of the right column is $s$ \er{whole}. The first square commutes by \er{lone}, and the rest are obviously commutative. 





\end{proof}

The back square in \er{cube2} is commutative by \er{whole}, and the other squares other than the front one are obviously commutative. Therefore the front square, which is just \er{delta2} is commutative. (Here we need to consider the composition of the front square with the skew morphism $C_*z(\ba^n)[1]\os{id}\to C_*z(\ba^n)[1]$.) Thus we have proved \er{ua}.
 
Therefore in diagram \er{chasing2}, 
$$\d^*(\tau_* cl_{\ba^n}^{-1})=\tau_* (\d^* cl_{\ba^n}^{-1})=\tau_* (l\cdot h_n^{-1})=0.$$ 
Now since we have known the type of $C_*z(\ba^{nl}-\ba^n)^{C_l}$ in Corollary \ref{pncp}, 
one has the top term in the middle column in \er{chasing2}
$$Hom(C_*z(\ba^{nl}-\ba^n)^{C_l},\bz/l(n)[2n])=0.$$

By the exactness of the middle column in \er{chasing2}, there is a unique 
$$g_n\in Hom(C_*z(\ba^{nl})^{C_l},\bz/l(n)[2n]),$$
which satisfies \er{fa}.

Now we want to prove that this $g_n$ satisfies \er{fu}.
This is equivalent to proving the commutativity of the second square in the following diagram 
$$\xymatrix{
C_*z(\ba^n)\ar[r]^\Delta\ar[d]^{cl_{\ba^n}^{-1}} & C_*z(\ba^{nl})^{C_l}\ar[d]^{g_n}\ar[r]^r & C_*z(\ba^{nl}-\ba^n)^{C_l}\ar[d]^{h_n^{-1}}\ar[r]^\delta & C_*z(\ba^n)[1]\ar[d]^{cl_{\ba^n}^{-1}[1]}\\
\bz(n)[2n]\ar[r]^\tau & \bz/l(n)[2n]\ar[r]^\beta & \bz(n)[2n+1]\ar[r]^{\cdot l} & \bz(n)[2n+1]
}$$

By \er{fa}, we see that the first square commutes. 
By the theory of triangulated categories, there exists a morphism $\phi:C_*z(\ba^{nl}-\ba^n)^{C_l}\to \bz(n)[2n+1]$ for the third vertical arrow which makes both the second and the third squares commute. 
If we use $h_n^{-1}$, 
the third square commutes by \er{ua}.
By Corollary \ref{pncp}, we see 
$$Hom(C_*z(\ba^{nl}-\ba^n)^{C_l},\bz(n)[2n+1])=\bz$$
 generated by $h_n^{-1}$.  Therefore $\phi$ must be $h_n^{-1}$ to make the third square commute. Hence $h_n^{-1}$ makes the second square commute, which proves \er{fu}.
\end{proof}

\begin{rmk}
One can directly check that $l(r^*(h_n^{-1}))=0$ in diagram \er{chasing2}, which is implied by \er{fu}. That is, we want to show the composition in the following diagram
$$\xymatrix{
 & & \bz(n)[2n+1]\\
C_*z(\ba^{nl})^{C_l}\ar[r]^(.4)r & C_*z(\ba^{nl}-\ba^n)^{C_l}\ar[r]^{h_n^{-1}} & \bz(n)[2n+1]\ar[u]^{\cdot l}
}$$
is zero.

The above diagram embeds in the following diagram 
$$\xymatrix{
C_*z(\ba^{nl})\ar[r]^r & C_*z(\ba^{nl}-\ba^n)\ar[r]^{s^{-1}} & \bz(n)[2n+1]\\
C_*z(\ba^{nl})^{C_l}\ar[r]^r\ar[u]^{C_*j} & C_*z(\ba^{nl}-\ba^n)^{C_l}\ar[r]^{h_n^{-1}}\ar[u]^{C_*j} & \bz(n)[2n+1]\ar[u]^{\cdot l}
,}$$
where the first two vertical maps are induced by the inclusions of the fixed point set sheaves to the whole sheaves. The commutativity of the last square is by Lemma \ref{leftsq}. 

The upper-right composition is zero by dimensional reason since it factorizes through $C_*z(\ba^{nl})$, which is $\bz(nl)[2nl]$. 
\end{rmk}

\begin{thm}\label{pf} The  morphism constructed above in \er{definef}
\bean
g:C_*z(\ba^{nl})^{C_l}\io\bigoplus_{i=n}^{nl-1}\bz/l(i)[2i]\oplus\bz(nl)[2nl].\label{f}
\eean
is an isomorphism. 
\end{thm}

\begin{proof} We denote the right hand side of \er{f} by $G$, and we consider the following diagram of distinguished triangles
$$\xymatrix{
C_*z(\ba^n)\ar[r]^\Delta\ar[d]^{cl_{\ba^n}^{-1}} & C_*z(\ba^{nl})^{C_l}\ar[r]^r \ar[d]^g & C_*z((\ba^{nl}-\ba^n)/C_l)\ar[r]^\delta\ar[d]^{h^{-1}} & C_*z(\ba^n)[1]\ar[d]^{cl_{\ba^n}^{-1}[1]}\\
\bz(n)[2n]\ar[r] & G\ar[r] & H\ar[r] & \bz(n)[2n+1]
,}$$
where the morphisms in the bottom row are the obvious ones (inclusions and projections) except the following 
$$\bz(n)[2n]\os\tau\to \bz/l(n)[2n]\os{\beta}\to \bz(n)[2n+1]\os{\cdot l}\to \bz(n)[2n+1].$$

The first square is commutative by \er{fa}, the second by \er{fu} and definition \er{g>n}, and the third by \er{ua}. 

Since both $cl_{\ba^n}^{-1}$ \er{cla} and $h^{-1}$ (Corollary \ref{pncp}) are isomorphisms, we conclude by the five lemma.
\end{proof}

Now let's work with $\bz/l$ coefficients. Let
\bea
C:=\bz(i)[j]\os{\cdot l}\to \bz(i)[j],
\l{dcplx}
\eea
with the second term placed at position 0, be a (double) complex in $DM$, and hence its total complex is an object in $DM$. The following natural map 
$$\lambda:C\to \bz/l(i)[j]$$
given by 
\bea 
\xymatrix{\bz(i)[j]\ar[r]^{\cdot l} & \bz(i)[j]\ar[d]^\tau\\
 & \bz/l(i)[j],
}
\eea
where $\tau:\bz(i)[j]\to \bz/l(i)[j]$ is the coefficient reduction, is an isomorphism in $DM$. Therefore $C$ is a canonical model of $\bz/l(i)[j]$. 

\begin{lem}\l{onebeta}  
One has a canonical isomorphism 
$$\kappa: C\otimes \bz/l\to \bz/l(i)[j]\op \bz/l(i)[j+1].$$
Furthermore, one has the following commutative diagram 
$$\xymatrix{
C\ar[r]^(.27)\lambda\ar[d]^\tau & \bz/l(i)[j]\ar[d]^{(id,\beta)}\\
C\otimes \bz/l\ar[r]^(.27)\kappa & \bz/l(i)[j]\op \bz/l(i)[j+1]
,}$$
where $\beta:\bz/l(i)[j]\to \bz/l(i)[j+1]$ is the Bockstein associated to the exact sequence of sheaves 
$$0\to \bz/l\os{\cdot l}\to \bz/l^2\to \bz/l\to 0.$$
\end{lem}

\begin{proof} One has 
\bea\l{silly}
C\otimes \bz/l=\bz/l(i)[j]\os{\cdot l}\to \bz/l(i)[j].
\eea
The map $\cdot l$ is zero here. The natural maps
$$\xymatrix{\bz/l(i)[j]\ar[r]^{\cdot l} & \bz/l(i)[j]\ar[d]^=\\
 & \bz/l(i)[j]}, \xymatrix {\bz/l(i)[j]\ar[d]^=\ar[r]^{\cdot l} & \bz/l(i)[j]\\
  \bz/l(i)[j] & }$$
give the isomorphism 
$$\kappa:C\otimes \bz/l\os\cong \to \bz/l(i)[j]\op\bz/l(i)[j+1].$$ 
The natural coefficient reduction map 
$$\tau:C\to C\otimes \bz/l$$
is
$$\xymatrix{
\bz(i)[j]\ar[r]^{\cdot l}\ar[d]^\tau & \bz(i)[j]\ar[d]^\tau\\
\bz/l(i)[j]\ar[r]^{\cdot l} & \bz/l(i)[j]}.$$
The composition $\kappa\circ \tau$ is the direct sum of the following 
$$\xymatrix {\bz(i)[j]\ar[r]^{\cdot l} & \bz(i)[j]\ar[d]^\tau\\
 & \bz/l(i)[j] }, 
\xymatrix {\bz(i)[j]\ar[d]^\tau \ar[r]^{\cdot l} & \bz(i)[j],\\
 \bz/l(i)[j] & }$$ 
which by definition are $\lambda$ and $\beta\circ \lambda$. 
\end{proof}


\begin{cor} Notation as in Theorem \ref{pz}. The natural morphism
$$f_m^{\bz/l}: F_m\otimes \bz/l\to M^c(L_m,\bz/l),$$
which by Lemma \ref{onebeta} is
$$f_m^{\bz/l}: \bigoplus_{i=0}^{m-1} \bz/l(i)[2i+1]\op \bigoplus_{i=1}^{m} \bz/l(i)[2i]\to M^c(L_m,\bz/l),$$
is an isomorphism. 
\end{cor}

\begin{prop}\label{pf/l} Notation as Theorem \ref{pf}. The natural morphism 
$$g^{\bz/l}: C_*z^{\bz/l}(\ba^{nl})^{C_l}\to G\otimes \bz/l,$$
which by Lemma \ref{onebeta} is
\bean
g^{\bz/l}: C_*z^{\bz/l}(\ba^{nl})^{C_l}\to \bigoplus_{i=n}^{nl-1}(\bz/l(i)[2i]\oplus \bz/l(i)[2i+1])\op \bz/l(nl)[2nl]
\label{f/l}
\eean
is an isomorphism.

Furthermore, one has the following commutative diagram 
\bea\xymatrix{
C_*z(\ba^{nl})^{C_l}\ar[d]^\tau\ar[r]^(.25)g & \bigoplus_{i=n}^{nl-1}\bz/l(i)[2i]\oplus \bz(nl)[2nl]\ar[d]^{\bigoplus_{i=n}^{nl-1}(id,\beta)\oplus \tau}\\
C_*z^{\bz/l}(\ba^{nl})^{C_l}\ar[r]^(.25){g^{\bz/l}} & \bigoplus_{i=n}^{nl-1}(\bz/l(i)[2i]\oplus \bz/l(i)[2i+1])\op \bz/l(nl)[2nl],
}\eea
where the $\tau$'s are coefficient reductions. 

\end{prop}

Now we give an interpretation of the above proposition in the $\ba^1$-homotopy category $H^{\ba^1}_\bullet (k)$. Recall that in the Dold-Kan correspondence, we have a pair of adjoint equivalence functors $N$ and $K$ between the category of simplicial abelian sheaves and the category of chain complexes of abelian sheaves concentrated in positive dimensions. For a general complex, we need to first apply the good truncation functor at position 0. Let 
$$K(\bz/l,j,i)=K(\bz/l(i)[j])$$ 
be the simplicial abelian sheaf corresponding to the motivic complex $\bz/l(i)[j]$. Then considered as a simplicial sheaf of pointed sets, it serves as the motivic Eilenberg-MacLane space in $H^{\ba^1}_\bullet (k)$ (see \er{rep}). 

The functor $K$ respects finite direct sums, i.e. it takes finite direct sums of chain complexes to finite direct sums of simplicial abelian sheaves. The forgetful functor from the category of simplicial abelian sheaves to the category of sheaves of pointed sets takes finite direct sums to products. Therefore using Proposition \ref{pf/l}, we have a natural isomorphism
\bea
&&K(C_* z^{\bz/l}(\ba^{nl})^{C_l})\os\cong\to K(\bigoplus_{i=n}^{nl}\bz/l(i)[2i]\oplus\bigoplus_{i=n}^{nl-1}\bz/l(i)[2i+1])\\
&=& \prod_{i=n}^{nl}K(\bz/l,2i,i)\times\prod_{i=n}^{nl-1}K(\bz/l,2i+1,i).
\eea

By \cite[Lemma 2.5.2]{dv}, we have an $\ba^1$-weak equivalence from $K(z^{\bz/l}(\ba^{nl})^{C_l})$ to\\ $K(C_* z^{\bz/l}(\ba^{nl})^{C_l})$, where $z^{\bz/l}(\ba^{nl})^{C_l}$ is considered to be a complex concentrated at dimension 0. It is easy to see by the construction of $K$ that $K(z^{\bz/l}(\ba^{nl})^{C_l})$ is the constant simplicial sheaf at $z^{\bz/l}(\ba^{nl})^{C_l}$. Thus we arrive at the following theorem. 

\begin{thm}\label{phomtype} In $H^{\ba^1}_\bullet(k)$, one has a natural isomorphism
\bean
z^{\bz/l}(\ba^{nl})^{C_l}\cong \prod_{i=n}^{nl}K(\bz/l,2i,i)\times\prod_{i=n}^{nl-1}K(\bz/l,2i+1,i).\l{homtype}
\eean
\end{thm}
This implies that $z^{\bz/l}(\ba^{nl})^{C_l}$ has the desired homotopy type to be the target for the total reduced power operation (\ref{op}).  The above theorem is part 1 of our Main Theorem \ref{main}. 

Similarly, using Theorem \ref{pf} we also have the following proposition.
\begin{prop}\label{phomtypez} In $H^{\ba^1}_\bullet(k)$, we have a natural isomorphism
\bean
z(\ba^{nl})^{C_l}\cong \prod_{i=n}^{nl-1}K(\bz/l,2i,i)\times K(\bz,2nl,nl).\l{homtypez}
\eean
\end{prop}


\section{Bockstein homomorphisms}\l{3.4}

We have constructed the following total reduced power operation (\ref{op})
$${\cal P}: z^{\bz/l}(\ba^n)\to z^{\bz/l}(\ba^{nl})^{C_l},$$
which in the $\ba^1$-homotopy category by Theorem \ref{phomtype} is canonically
\begin{eqnarray*}
&&{{\cal P}}:K(\bz/l,2n,n)\to \prod_{i=n}^{nl} K(\bz/l, 2i,i)\times \prod_{i=n}^{nl-1} K(\bz/l, 2i+1,i)\\
&=& \prod_{j=0}^{n(l-1)} K(\bz/l, 2n+2j,n+j)\times \prod_{j=0}^{n(l-1)-1} K(\bz/l, 2n+2j+1,n+j).
\eea

We introduce the following notation: let
\bean
D^i: K(\bz/l,2n,n)\to K(\bz/l,2n+2i,n+i),\ 0\leq i\leq n(l-1)\l{di}
\eean
and 
\bean
E^i: K(\bz/l,2n,n)\to K(\bz/l,2n+2i+1,n+i),\ 0\leq i\leq n(l-1)-1\l{ei}
\eean
denote the components of our total reduced power operation $\cp$.

The goal of this section is to prove Theorem \ref{pbock}, which is part 2 of our Main Theorem \ref{main}. 

First we want to prove that our operations with $\bz/l$ coefficients factorizes through a modified integral version. 

Consider the following map 
$$l\varpi: z(\ba^{nl})\to z(\ba^{nl})^{C_l}$$
defined as follows. Given a smooth scheme $U$ and a cycle $Z\in z(\ba^{nl})(U)$, define 
$$l\varpi(Z)=l\sum_{\sigma\in C_l} \sigma(Z),$$ 
where $\sigma(Z)$ is the image of $Z$ under an element $\sigma\in C_l$. It is clear that $l\varpi(Z)$ is invariant under the $C_l$ action by its construction and thus is in $z(\ba^{nl})^{C_l}(U)$. Denote by $z(\ba^{nl})^{C_l}/{\\\rm Im}(l\varpi)$ the Nisnevich sheaf associated to the presheaf whose value on a smooth scheme $U$ is $z(\ba^{nl})^{C_l}(U)/{\rm Im}(l\varpi:z(\ba^{nl})(U)\to z(\ba^{nl})^{C_l}(U))$, where the quotient is in the category of abelian groups. 


\begin{lem}\l{proj} In the $\ba^1$-homotopy category $H^{\ba^1}_\bullet (k)$, one has a natural isomorphism
\bean
q:z(\ba^{nl})^{C_l}/{\rm Im}(l\varpi)\os\cong\to \prod_{i=n}^{nl-1} K(\bz/l,2i,i)\times K(\bz/l^2,2nl,nl).\l{zq}
\eean

Furthermore, one has the following commutative diagram 
\bean
\xymatrix{
z(\ba^{nl})^{C_l}\ar[d]^p\ar[r]^(.3)g & \prod_{i=n}^{nl-1} K(\bz/l,2i,i)\times K(\bz,2nl,nl)\ar[d]^{\prod_{i=n}^{nl-1} id\times \eta}\\
z(\ba^{nl})^{C_l}/{\rm Im}(l\varpi)\ar[r]^(.3)q & \prod_{i=n}^{nl-1} K(\bz/l,2i,i)\times K(\bz/l^2,2nl,nl),
}\l{qim}
\eean
where $p$ is the natural projection, the two horizontal isomorphisms $g$ and $q$ are by \er{homtypez} and \er{zq}, and $\eta$ on the right vertical arrow is the natural coefficient reduction map. 
\end{lem}

\begin{proof} Apply the singular complex functor $C_*$ and consider the following statement on the motive level corresponding to \er{zq}:
\bean\l{q}
q:C_*(z(\ba^{nl})^{C_l}/{\rm Im}(l\varpi))\os\cong \to\oplus_{i=n}^{nl-1} \bz/l(i)[2i]\oplus \bz/l^2(nl)[2nl].
\eean

One has the following natural distinguished triangle
\bean
C_* z(\ba^{nl})\os {C_*(l\varpi)}\longrightarrow C_*z(\ba^{nl})^{C_l}\os {C_*p} \to C_*(z(\ba^{nl})^{C_l}/{\rm Im}(l\varpi))\to C_* z(\ba^{nl})[1].\l{lomg}
\eean

To define the $i$-th component 
$$q_i: C_*(z(\ba^{nl})^{C_l}/{\rm Im}(l\varpi))\to \bz/l(i)[2i],\ n\leq i\leq nl-1$$
of $q$, we study the long exact sequence of $Hom$ group associated to \er{lomg}
\bea
&0=&Hom(C_* z(\ba^{nl})[1],\bz/l(i)[2i])\to Hom(C_*(z(\ba^{nl})^{C_l}/{\rm Im}(l\varpi)),\bz/l(i)[2i])\\
&\to& Hom(C_*z(\ba^{nl})^{C_l},\bz/l(i)[2i])\to Hom(C_* z(\ba^{nl}),\bz/l(i)[2i])=0.
\eea
We pick $q_i\in Hom(C_*(z(\ba^{nl})^{C_l}/{\rm Im}(l\varpi)),\bz/l(i)[2i])$ to be the unique element corresponding to $g_i\in Hom(C_*z(\ba^{nl})^{C_l},\bz/l(i)[2i])$, which is the $i$-th component of $g$ in \er{f}.  

The map $C_*(l\varpi): C_*z(\ba^{nl})\to C_*z(\ba^{nl})^{C_l}$ can be seen to induce the multiplication by $l^2$ on the $\bz(nl)[2nl]$ factor for the natural motive types (see \er{cla}, \er{f}). Using a diagram similar to \er{chasing2}, one can define a unique morphism
$$q_{nl}: C_*(z(\ba^{nl})^{C_l}/{\rm Im}(l\varpi))\to \bz/l^2(nl)[2nl],$$
such that the following diagram is commutative
$$\xymatrix{
C_*z(\ba^{nl})^{C_l}\ar[r]^(.4)p\ar[d]^{g_{nl}} & C_*(z(\ba^{nl})^{C_l}/{\rm Im}(l\varpi))\ar[d]^{q_{nl}}\\
\bz(nl)[2nl]\ar[r]^(.4)\eta & \bz/l^2(nl)[2nl].
}$$

Then one can conclude that the above constructed $q$ is an isomorphism by a diagram of distinguished triangles and the five lemma similar to the proof of Theorem \ref{pf}. The commutativity of the motive version of \er{qim} is clear from our construction of $q$. 

The same analysis as in the three paragraphs before Theorem \ref{phomtype}, using the functor $K$ in the Dold-Kan correspondence and \cite[Lemma 2.5.2]{dv}, gives our statements on the homotopy category level. 
\end{proof}

We now want to prove that our total reduced power operation $\cal P$ \er{op} for $\bz/l$ coefficients factors through $z(\ba^{nl})^{C_l}/{\rm Im}(l\varpi)$ with $\bz$ coefficients. Namely we will show the existence of the following commutative diagram
\bean
\xymatrix{
z(\ba^{nl})^{C_l}/{\rm Im}(l\varpi)\ar[dr]^{\tau'} & z(\ba^{nl})^{C_l}\ar[l]_p\ar[d]^\tau\\
z^{\bz/l}(\ba^n)\ar[r]^{{\cal P}}\ar[u]^{\wp} & z^{\bz/l}(\ba^{nl})^{C_l}.\l{dwp}
}
\eean
where $\tau$ is the coefficient reduction by $l$, $\tau'$ is defined because $Im(l\varpi)$ is reduced to zero. 

The proof of the following proposition is similar to that of \cite[Proposition 8.1]v.

\begin{prop}\l{lift} The lift $\wp$ in \er{dwp} exists. 
\end{prop}

\begin{proof} Recall first that our total reduced power operation ${\cal P}$ \er{op} is also defined for $\bz$ coefficients, and we have ${\cal P}: z(\ba^n)\to z(\ba^{nl})^{C_l}$. A cycle in $z^{\bz/l}(\ba^n)(U)$ lifts to an integral cycle in $z(\ba^n)(U)$ unique up to a multiple of $l$. To prove the existence of $\wp$, we only need to show that for $Z_1,Z_2\in z(\ba^n)(U)$ such that $Z_1-Z_2=lW$ for some $W\in  z(\ba^n)(U)$, we have ${\cal P}(Z_1)-{\cal P}(Z_2)\in Im(l\varpi)$, i.e. 
$$Z_1^l-Z_2^l=l\varpi(T)$$
for some $T\in z(\ba^{nl})(U)$. 
Observe that any invariant cycle in $z(\ba^{nl})^{C_l}(U)$ divisible by $l^2$, say $l^2 V$ with $V$ invariant, is equal to $l\varpi(V)$. 
Because $Z_1^l-Z_2^l$ is invariant, we only need to take care of the terms in $Z_1^l-Z_2^l=(Z_2+lW)^l-Z_2^l$ that are not divisible by $l^2$. Such terms are $l\sum_i Z_2\times\cdots\times \underset {\rm i-th\ slot} W\times\cdots\times Z_2$, and it can be written as $l\varpi(W\times Z_2\times\cdots\times Z_2)$. 
\end{proof}

Now we can prove our result concerning the Bockstein homomorphism. 

\begin{thm}\l{pbock} $E^i=\beta\circ D^i$ for $0\leq i\leq nl-1$.
\end{thm}

\begin{proof} We have a natural factorization (see \er{dwp})
$$\tau:z(\ba^{nl})^{C_l}\os p \to z(\ba^{nl})^{C_l}/{\rm Im}(l\varpi)\os {\tau'} \to z^{\bz/l}(\ba^{nl})^{C_l}$$
of the coefficient reduction map. Lemma \ref{proj} says that $p$ induces the identities on the factors $K(\bz/l,2i,i)$ for $n\leq i\leq nl-1$, and thus $\tau'$ is the same as $\tau$ on these factors. Interpreting Proposition \ref{pf/l} in the homotopy category, we see that the $i$-th component $\tau_i$ of $\tau$ for $0\leq i\leq nl-1$ is 
$$\tau_i=(id,\beta):K(\bz/l,2i,i)\to K(\bz/l,2i,i)\times K(\bz/l,2i+1,i),$$
and therefore this also holds for $\tau'$. 
Proposition \ref{lift} gives the following factorization 
$${\cal P}: z^{\bz/l}(\ba^n)\os {\wp} \to z(\ba^{nl})^{C_l}/{\rm Im}(l\varpi)\os {\tau'}\to z^{\bz/l}(\ba^{nl})^{C_l}.$$
This proves our theorem. 
\end{proof}

\section{Basic properties}\l{3.5}


In this section, we prove some very basic properties of our operations $D^i$ \er{di} and $E^i$ \er{ei}. By the representability of motivic cohomology \er{rep}, our operations represent the following cohomology operations:
$$D^i: H^{2n,n}(-,\bz/l)\to H^{2n+2i,n+i}(-,\bz/l),\ 0\leq i\leq n(l-1);$$
$$E^i: H^{2n,n}(-,\bz/l)\to H^{2n+2i+1,n+i}(-,\bz/l),\ 0\leq i\leq n(l-1)-1.$$

Since our constructions of the $D^i$ and $E^i$ are through maps between the corresponding motivic Eilenberg-MacLane spaces in the $\ba^1$-homotopy category, these operations are natural with respect to morphisms between schemes (actually any morphism in $H^{\ba^1}_\bullet(k)$). 

\begin{thm} \l{p0} $D^0:K(\bz/l,2n,n)\to K(\bz/l,2n,n)$ is the identity.  
\end{thm}

\begin{proof} By our construction (see the proof of Proposition \ref{ng}, especially \er{fa}), the $K(\bz/l,2n,n)$ factor of $z^{\bz/l}(\ba^{nl})^{C_l}$ is represented by $z^{\bz/l}(\ba^n)$, the sheaf of equidimensional cycles supported on the diagonal $\ba^n$ of $\ba^{nl}$. Suppose $Z=\sum n_i Z_i\in z^{\bz/l}(\ba^n)(U)$, where $U$ is a smooth scheme and the $Z_i$ are distinct irreducible subschemes. Then 
$${\cal P}(Z)=\sum n_{i_1}\cdots n_{i_l} Z_{i_1}\times_U\cdots \times_U Z_{i_l}\in z^{\bz/l}(\ba^{nl})^{C_l}(U).$$
For $i_1,\cdots,i_l$ not all equal, $Z_{i_1}\times_U\cdots \times_UZ_{i_l}$ must have non-empty intersection with\\ 
$U\times (\ba^{nl}-\ba^n)$ (we assume that the $Z_i$ are  all distinct). For each $i$, the fiber product\\ $Z_{i}\times_U\cdots \times_U Z_{i}$ may be reducible. The following subvariety, denoted by $W_i$, of  $Z_{i}\times_U\cdots \times_U Z_{i}$ defined by $y_1=\cdots=y_l$, where the $y_j$ are the coordinates to the $j$-th $\ba^n$ component of $\ba^{nl}$, is obviously isomorphic to $Z_i$ under the identification of the diagonal of $\ba^{nl}$ to $\ba^n$, and so irreducible. By definition all other irreducible components of $Z_{i}\times_U\cdots \times_U Z_{i}$ have non-empty intersections with $U\times (\ba^{nl}-\ba^n)$. The coefficient of the irreducible subvariety $W_i$ is the same as the coefficient of $Z_{i}\times_U\cdots \times_U Z_{i}$, which is $n_i^l$. It is well known that $n_i^l=n_i$ in $\bz/l$. Therefore the part of $\cp(Z)$ supported on the diagonal of $\ba^{nl}$ is the same as the cycle $Z$. This completes the proof.
\end{proof}

\begin{thm} \l{lastp} $D^{n(l-1)}=P:K(2n,n,\bz/l)\to K(2nl,nl,\bz/l)$, where $P$ is the $l$-th cup product power operation \er{P}. 
\end{thm}

\begin{proof} By definition \er{factorization}, 
$$P=j\circ \cp:z^{\bz/l}(\ba^n)\to  z^{\bz/l}(\ba^{nl})^{C_l}\to z^{\bz/l}(\ba^{nl}).$$ 
Here 
$$j: z^{\bz/l}(\ba^{nl})^{C_l}\to z^{\bz/l}(\ba^{nl})$$ 
is the natural inclusion of the fixed point set sheaf into the whole sheaf, which in the $\ba^1$-homotopy category by Theorem \ref{phomtype} is canonically
$$j:\prod_{i=n}^{nl}K(\bz/l,2i,i)\times\prod_{i=n}^{nl-1}K(\bz/l,2i+1,i)\to K(\bz/l,2nl,nl).$$
We now prove that $j$ is isomorphic to the projection to the factor $K(\bz/l,2nl,nl)$. 


Let's consider the corresponding map on the motive level 
$$C_*j: C_*z^{\bz/l}(\ba^{nl})^{C_l}\to C_*z^{\bz/l}(\ba^{nl}).$$ 
In view of Proposition \ref{pf/l}, the map $C_*j$ canonically has the following form 
$$C_*j:\oplus_{i=n}^{nl-1}(\bz/l(i)[2i]\op \bz/l(i)[2i+1])\oplus \bz(nl)[2nl]\to \bz(nl)[2nl].$$ 
$C_*j$ can be seen to be zero on the factors $(\bz/l(i)[2i]\op \bz/l(i)[2i+1])$ for $n\leq i\leq nl-1$ by our construction. 
We now concentrate on the top dimensional factor $\bz/l(nl)[2nl]$. Consider the following diagram 
\bea
\xymatrix{
C_*z^{\bz/l}(\ba^{nl})^{C_l}\ar[r]^r\ar[d]^{C_*j} & C_*z^{\bz/l}(\ba^{nl}-\ba^n)^{C_l}\ar[r]^{\widetilde{\pi_*}}\ar[d]^{C_*j} & C_*z^{\bz/l}((\ba^{nl}-\ba^n)/C_l)\ar[d]^{\pi^*} \\
C_*z^{\bz/l}(\ba^{nl})\ar[r]^r & C_*z^{\bz/l}(\ba^{nl}-\ba^n)\ar[r]^=& C_*z^{\bz/l}(\ba^{nl}-\ba^n)\\
}
\eea
where $\widetilde{\pi_*}$ is that in \er{free}, and all the other arrows are the obvious ones. 
Observe that all the horizontal arrows induce the identities on the factor $\bz/l(nl)[2nl]$. Now we consider the effect of the last vertical arrow. The $\bz/l(nl)[2nl]$ factors of both  $C_*z^{\bz/l}((\ba^{nl}-\ba^n)/C_l)$ and $C_*z^{\bz/l}(\ba^{nl}-\ba^n)$ are represented by the fundamental classes. We see that $\pi^*$ 
is the identity on this factor by Lemma \ref{l-rest}. Therefore so is the first vertical one. 

\end{proof}

\section{Comparison with Voevodsky's operations}\l{4.2}

We start by explaining the following commutative diagram
\bean
\xymatrix{
z({\b A}^n)\ar[r]^{{\cal P}} \ar[dr]^{\Psi} & z({\b A}^{nl})^{C_l}\ar[d]^{\gamma}\\
  & z({\b A}^{nl})^{hC_l}.
}\l{comparison}
\eean
where ${\cal P}: z(\ba^n)\to z(\ba^{nl})^{C_l}$ is our total reduced power operation \er{op}. 

$z({\b A}^{nl})^{hC_l}$ in (\ref{comparison}) is an analogue of the homotopy fixed point set in topology. To define it, first note that for any linear algebraic group $G$, Voevodsky has defined its classifying space $BG$ in the $\ba^1$-homotopy category $H^{\ba^1}_\bullet(k)$ (see \cite[Section 6]{v} and \cite{mv}). For the reader's convenience, let's recall its definition. 

Suppose that $G\to GL(V)$ is a faithful representation of $G$. Denote by $\widetilde{V_i}$ the open subset in $\ba(V)^i$ where $G$ acts freely. We have a sequence of closed embeddings $f_i:\widetilde{V_i}\to \widetilde{V_{i+1}}$ given by $(v_1,\cdots,v_i)\to (v_1,\cdots,v_i,0)$. Set $EG={\rm colim}_i \tvi$ and $BG={\rm colim}_i \widetilde{V_i}/G$, where $\widetilde{V_i}/G$ is the quotient scheme, and the colimits are taken in the category of sheaves. It is known that the homotopy type of $BG$ doesn't depend on the choice of the representation $G\to GL(V)$. In particular, we have the classifying space $BC_l$ for the cyclic group. Under our assumption of the existence of a primitive $l$-th root of unity $\zeta$ in $k$, we have an isomorphism between $C_l$ and $\mu_l$, the group of $l$-th roots of unity in $k$. Therefore we have an $\ba^1$-weak equivalence $BC_l\simeq B\mu_l$. 

Now let's define $z({\b A}^{nl})^{hC_l}$. Fix a representation $C_l\to GL(V)$, and use the same notation as above.

\begin{defn}\l{defn} Given a smooth scheme Y and a scheme X, define $z(Y,X)$ to be the sheaf whose value on a smooth scheme $U$ is $z(Y,X)(U):=z(X)(Y\times U)$. 

Define $z(\widetilde{V_i},\ba^{nl})^{C_l}$ to be the sheaf whose value $z(\widetilde{V_i},\ba^{nl})^{C_l}(U)$ on a smooth scheme $U$ is the group of cycles in $z(\ba^{nl})(\widetilde{V_i}\times U)$ which are invariant under the $C_l$ action on $\ba^{nl}\times \widetilde{V_i}\times U$ as the product of its natural actions on $\ba^{nl}$ and $\widetilde{V_i}$ and the trivial action on $U$. Concretely, a cycle $Z$ belongs to $z(\widetilde{V_i},\ba^{nl})^{C_l}(U)$ if for any $\sigma\in C_l$, $v\in \tvi$ and $u\in U$, the fiber of $Z$ over the point $(\sigma(v),u)\in \tvi\times U$, as a linear combination of points in $\ba^{nl}$, is the image of the fiber of $Z$ over $(v,u)$ under $\sigma$. 

The natural inclusion $f_i:\widetilde{V_i}\to \widetilde{V_{i+1}}$ induces a natural pullback map 
$$Cycl(f_i):z(\widetilde{V_i},\ba^{nl})\gets z(\widetilde{V_{i+1}},\ba^{nl}),$$ 
and it can be seen to induce a map on the fixed point set 
\bean
Cycl(f_i):z(\widetilde{V_i},\ba^{nl})^{C_l}\gets z(\widetilde{V_{i+1}},\ba^{nl})^{C_l},
\l{cyclfi}
\eean
since $\widetilde{V_i}\times 0\subset \widetilde{V_{i+1}}$ is invariant under the $C_l$ action. 

We define 
$$z(\ba^{nl})^{hC_l}:={\rm lim}_i z(\widetilde {V_i},\ba^{nl})^{C_l},$$ 
where the limit is taken in the category of sheaves.

Furthermore we define $z^{\bz/l}(\ba^{nl})^{hC_l}$ to be $z(\ba^{nl})^{hC_l}\otimes \bz/l$. 

\end{defn}

Since the action of $C_l$ on $\widetilde{V_i}$ is free, so is its action on the product $\ba^{nl}\times \tvi\times U$ for any scheme $U$. Then the projection 
$$p:\ba^{nl}\times_{C_l} \tvi\times U:=(\ba^{nl}\times \tvi\times U)/C_l\to (\tvi\times U)/C_l=\tvi/C_l\times U$$ 
is a vector bundle of rank $nl$. Furthermore any cycle $Z\in z(\widetilde{V_i},\ba^{nl})^{C_l}(U)$, being invariant under the free $C_l$ action on $\ba^{nl}\times \tvi\times U$, comes from a cycle $Z'$ on the quotient (i.e., the flat pullback of $Z'$ is $Z$, see \er{free}). As a cycle on the vector bundle 
$$\ba^{nl}\times_{C_l} \tvi\times U\os p \to \tvi/C_l\times U,$$ 
$Z'$ is equidimensional of relative dimension 0 over the base $\tvi/C_l\times U$ under the projection $p$, since originally $Z$ is equidimensional over $\tvi\times U$. 

\begin{rmk}\l{informal} One has the following imprecise interpretation of $z(\ba^{nl})^{hC_l}$. Given a smooth scheme $U$, $z({\b A}^{nl})^{hC_l}(U)$ can be thought of as the free abelian group on the  closed irreducible subschemes on ${\b A}^{nl}\times_{C_l} EC_l\times U$, which are equidimensional of relative dimension 0 
over $BC_l\times U$ under the vector bundle projection $p: {\b A}^{nl}\times_{C_l} EC_l\times U \to BC_l\times U$. This is imprecise because $EC_l$ and $BC_l$ are not schemes, but rather sheaves. 
\end{rmk}


$\gamma: z({\b A}^{nl})^{C_l}\to z({\b A}^{nl})^{hC_l}$ in (\ref{comparison}) is an analogue of the inclusion of a fixed point set into its associated homotopy fixed point set in topology, and is defined as follows. For a cycle $Z\in z({\b A}^{nl})^{C_l}(U)$, $Cycl(pr_i)(Z)\in z({\b A}^{nl})(\tvi\times U)$ is the pullback of $Z$ by the projection $pr_i: \tvi\times U\to U$. It is easily seen that $Cycl(pr_i)(Z)\in z(\tvi,{\b A}^{nl})^{C_l}(U)$, since $Z\in z({\b A}^{nl})^{C_l}(U)$. Furthermore 
$$Cycl(f_i)(Cycl(pr_{i+1})(Z))=Cycl(pr_i)(Z)$$ 
with $Cycl(f_i)$ as in \er{cyclfi}, 
since clearly 
$$pr_{i+1}\circ (f_i\times id_U)=pr_i:\tvi\times U\to U.$$ 
Therefore the system $\{Cycl(pr_i)(Z)\}_i$ defines a cycle in the limit $z(\ba^{nl})^{hC_l}(U)$, and this is defined to be $\gamma (Z)$. 

In \er{comparison}, $\Psi=\gamma\circ {\cal P}$ is the composition. It is understood that $\Psi$ is the same as Voevodsky's construction \cite[Section 5]{v}. 

Now let's concentrate on $\bz/l$ coefficients. It is a computation of Voevodsky \cite[Sections 5 and 6] v that one has a canonical isomorphism
\bean
z^{\bz/l}(\ba^{nl})^{hC_l}\cong \prod_{i=0}^{nl} K(\bz/l,2i,i)\times \prod_{i=0}^{nl-1} K(\bz/l,2i+1,i).\l{hcl}
\eean
To be more precise, he (loc. cit.) first proved that for a smooth scheme $X$ 
\bean
Hom_{H^{\ba^1}_\bullet}(X_+,z^{\bz/l}({\b A}^{nl})^{hC_l})\cong H^{2nl,nl}(X\times BC_l,\bz/l).\l{bdlclass}
\eean

\begin{rmk} \er{bdlclass} can be informally understood as follows. 
In view of Remark \ref{informal}, a cycle $Z\in z^{\bz/l}({\b A}^{nl})^{hC_l}(X)$ "should" define a motivic cohomology class of dimension $(2nl,nl)$ of the base $X\times BC_l$, since if the vector bundle $\ba^{nl}\times_{C_l} EC_l\times X$ were trivial, i.e. $\ba^{nl}\times BC_l\times X$, such a cycle does define such a class by definition (see \er{rep} and \er{em}). Actually much of Voevodsky's proof of \er{bdlclass} is to add another bundle to make the direct sum trivial, and then use the Thom isomorphism to "pull" the dimension of the motivic cohomology class back. 
\end{rmk}

Under our assumption, 
we have a weak equivalence $BC_l\simeq B\mu_l$.
Then Voevodsky computed the motivic cohomology of $X\times B\mu_l$ with $\bz/l$ coefficients. His result (loc. cit.) is 
\bean
H^{2nl,nl}(X\times B\mu_l,\bz/l)\cong \bigoplus_{\epsilon=0,1,\ 0\leq \delta\leq nl} H^{*,*}(X,\bz/l)u^\epsilon v^\delta,\l{h2nl}
\eean
as an $H^{*,*}(X,\bz/l)$ module. Here $v\in H^{2,1}(B\mu_l,\bz/l)$ corresponds to the line bundle on $B\mu_l$ associated to the standard representation $\rho$, and $u\in H^{1,1}(B\mu_l,\bz/l)$ is the unique class such that its Bockstein $\beta(u)=v$ and its restriction to a rational point $*$ of $B\mu_l$ which lifts to a rational point of some $\tvi$ is zero. In the following table, we count the possible dimensions of $u^\epsilon v^\delta$ for $\epsilon=0,1,\ 0\leq \delta\leq nl$ and calculate the corresponding dimensions of $H^{*,*}(X,\bz/l)$ such that the sum of the dimensions is $(2nl,nl)$: 
\bean
&&\begin{tabular}{|r|r|r|r|}\hline 
$\delta$&0&0&1\\\hline
$\epsilon$ &0 &1 &0\\\hline
${\rm dim}u^\epsilon v^\delta$ &(0,0) &(1,1) &(2,1) \\\hline
${\rm dim}H^{*,*}(X)$ &$(2nl,nl)$ & $(2nl-1,nl-1)$&$(2nl-2,nl-1)$\\\hline
\end{tabular}\nm\\
&&\begin{tabular}{|r|r|r|r|}\hline 
1&$\cdots$&$nl$&$\delta$\\\hline
1&$\cdots$&0&$\epsilon$\\\hline
(3,2)&$\cdots$&$(2nl,nl)$&${\rm dim}u^\epsilon v^\delta$ \\\hline
$(2nl-3,nl-2)$&$\cdots$&(0,0)&${\rm dim}H^{*,*}(X)$\\\hline
\end{tabular}\nm
\eean

Combining \er{bdlclass}, \er{h2nl} and the above table of possible dimensions of $H^{*,*}(X,\bz/l)$, we arrive at the conclusion \er{hcl}.  

Recall our Theorem \ref{phomtype} about the homotopy type of $z^{\bz/l}(\ba^{nl})^{C_l}$. To compare the two operations, we need to study the map $\gamma: z^{\bz/l}(\ba^{nl})^{C_l}\to z^{\bz/l}(\ba^{nl})^{hC_l}$ in \er{comparison}. Our result is 

\begin{thm}\l{tcomp} The map $\gamma: z^{\bz/l}(\ba^{nl})^{C_l}\to z^{\bz/l}(\ba^{nl})^{hC_l}$ in \er{comparison}, which in view of Theorem \ref{phomtype} and \er{hcl} has the form
\bean
&\gamma:& \prod_{i=n}^{nl}K(\bz/l,2i,i)\times\prod_{i=n}^{nl-1}K(\bz/l,2i+1,i)\nm\\
&\to& \prod_{i=0}^{nl} K(\bz/l,2i,i)\times \prod_{i=0}^{nl-1} K(\bz/l,2i+1,i),\nm
\eean
induces isomorphisms on the factors $K(\bz/l,2i,i)$ for $n\leq i\leq nl$. In view of the diagram (\ref{comparison}), we see that our operations $D^i$ \er{di} coincide with the corresponding ones of Voevodsky for cohomology classes of dimension $(2n,n)$, up to nonzero constants in $\bz/l$. 
\end{thm}

\begin{rmk} In view of Theorem \ref{pbock} and \cite[Lemma 9.6]v, we see that our $E^i$ \er{ei} are also the same, up to nonzero constants, as the corresponding ones of Voevodsky for cohomology classes of dimension $(2n,n)$. 
\end{rmk}

Actually the above theorem will be a special case of the following general situation. We now adopt the following notational convention in the rest of this section. We assume that $\ba^n$ always has a trivial $\mu_l$ action, $\ba^m$ has a $\mu_l$ action which is a direct sum of nontrivial irreducible representations, and $\ba^{n+m}$ is the direct sum of them. 

Our previous computations in Sections \ref{3.2} and \ref{3.3}, which correspond to the situation that $m=n(l-1)$, generalize to this general situation. An analogue of Theorem \ref{phomtype} gives us a natural isomorphism
\bean
z^{\bz/l}(\ba^{n+m})^{\mu_l}\cong \prod_{i=n}^{n+m} K(\bz/l,2i,i)\times \prod_{i=n}^{n+m-1} K(\bz/l,2i+1,i).\l{mul}
\eean

Denote by $z^{\bz/l}(\ba^{n+m})^{h\mu_l}$ the similarly defined homotopy fixed point set as in Definition \ref{defn}. The following general versions of \er{bdlclass}, \er{h2nl} and \er{hcl}:
\bean
Hom_{H^{\ba^1}_\bullet}(X_+,z^{\bz/l}({\b A}^{n+m})^{h\mu_l})\cong H^{2(n+m),n+m}(X\times B\mu_l,\bz/l),\l{bcnm}
\eean
\bean
H^{2(n+m),n+m}(X\times B\mu_l,\bz/l)\cong \bigoplus_{\epsilon=0,1,\ 0\leq \delta\leq n+m} H^{*,*}(X,\bz/l)u^\epsilon v^\delta,\l{hnm}
\eean
and 
\bean
z^{\bz/l}(\ba^{n+m})^{h\mu_l}\cong \prod_{i=0}^{n+m} K(\bz/l,2i,i)\times \prod_{i=0}^{n+m-1} K(\bz/l,2i+1,i)\l{hmul}
\eean
still hold.

Similar to diagram \er{comparison}, we can define a natural map $\gamma: z^{\bz/l}(\ba^{n+m})^{\mu_l}\to z^{\bz/l}(\ba^{n+m})^{h\mu_l}$. We now prove the following theorem, which in particular gives Theorem \ref{tcomp}. 

\begin{thm}\l{mn} In view of \er{mul} and \er{hmul}, $\gamma: z^{\bz/l}(\ba^{n+m})^{\mu_l}\to z^{\bz/l}(\ba^{n+m})^{h\mu_l}$ induces isomorphisms on the factors $K(\bz/l,2i,i)$ for $n\leq i\leq m+n$. 
\end{thm}

First a lemma. 

\begin{lem}\l{ltop} For arbitrary dimensions $n$ and $m$, the map $\gamma$ induces the identity on the top dimensional factor $K(\bz/l,2(n+m),n+m)$.
\end{lem}

\begin{proof} Consider the following diagram 
\bean
\xymatrix{
z^{\bz/l}(\ba^{n+m})^{\mu_l}\ar[r]^j\ar[d]^{\gamma} & z^{\bz/l}(\ba^{n+m})\\
z^{\bz/l}(\ba^{n+m})^{h\mu_l}\ar[ur]^{j'}, & 
}\l{top}
\eean
where $j$ is the inclusion of the fixed point set sheaf into the whole sheaf, and $j'$ is defined as follows. Regard $z^{\bz/l}(\ba^{n+m})^{h\mu_l}$ as $z^{\bz/l}(E\mu_l,\ba^{n+m})^{\mu_l}$ (see Definition \ref{defn}). Let $incl:*\to E\mu_l$ be a rational point of $E\mu_l$ which lifts to a rational point of some $\tvi$. Then define 
$$j':z^{\bz/l}(E\mu_l,\ba^{n+m})^{\mu_l}\to z^{\bz/l}(E\mu_l,\ba^{n+m})\os {Cycl(incl)} \longrightarrow z^{\bz/l}(\ba^{n+m})$$ 
to be the composition of the inclusion map and the pullback map by $incl$. Note that $j'$ doesn't depend on the choice of the rational point $*$ since $E\mu_l$ is contractible. 

The diagram \er{top} can be seen to commute as follows. Observe that $\gamma=Cycl(pr)$ is the pullback map by $pr:E\mu_l\to *$. Therefore 
$$j'\circ \gamma=Cycl(incl)\circ Cycl(pr)=Cycl(pr\circ incl)=Cycl(id_*)$$ 
is the identity map, except that now we forget the fixed point set structure, which is exactly $j$. 

Exactly the same arguments as in the proof of Theorem \ref{lastp} show that $j$ induces the identity on the factor $K(\bz/l,2(n+m),n+m)$. We now want to prove that $j'$ also induces the identity on this top dimensional factor.  For a scheme $X$, we have the following commutative diagram:
$$\xymatrix{
Hom_{H^{\ba^1}_\bullet}(X_+,z^{\bz/l}({\b A}^{n+m})^{h\mu_l})\ar[r]^\cong \ar[d]^{j'\circ} & H^{2(n+m),n+m}(X\times B\mu_l,\bz/l)\ar[d]^{j'_H}\\
Hom_{H^{\ba^1}_\bullet}(X_+,z^{\bz/l}({\b A}^{n+m}))\ar[r]^\cong  & H^{2(n+m),n+m}(X,\bz/l),
}$$
where the first isomorphism is \er{bcnm}, and the second isomorphism is the obvious one (see \er{rep} and \er{em}). The left vertical map $j'\circ$ is composition by $j'$, and the right vertical map $j'_H$ is induced by the inclusion $*\to E\mu_l\to B\mu_l$. The diagram commutes by the naturality of \er{bcnm} (see \cite[Lemma 5.12]v). \er{hnm} says that $j'_H$ has the form 
\bean
j'_H:H^{2(n+m),n+m}(X\times B\mu_l,\bz/l)&\cong& \bigoplus_{\epsilon=0,1,\ 0\leq \delta\leq n+m} H^{*,*}(X,\bz/l)u^\epsilon v^\delta\nm\\
&\to& H^{2(n+m),n+m}(X,\bz/l).\nm
\eean
Actually $j'_H$ is the projection to the summand $H^{2(n+m),n+m}(X,\bz/l)$, which corresponds to $\epsilon=\delta=0$ in $\bigoplus_{\epsilon=0,1,\ 0\leq \delta\leq nl} H^{*,*}(X,\bz/l)u^\epsilon v^\delta$, since both $v$ (by dimension reasons) and $u$ (by definition) are restricted to zero in the motivic cohomology of the rational point $*$. This means that $j':z^{\bz/l}(\ba^{n+m})^{h\mu_l}\to z^{\bz/l}(\ba^{n+m})$ is actually the projection to the factor $K(\bz/l,2(n+m),n+m)$. 

By the commutativity of the diagram \er{top}, we see that $\gamma$ is the identity on the factor $K(\bz/l,2(n+m),n+m)$. 
\end{proof}

Now let's prove our theorem. 


\begin{pfmn} We run induction on $m$. When $m=0$ we are done by Lemma \ref{ltop}. 
Now assume that the assertion holds for $m-1$, i.e. (with obvious notation) 
$$\gamma_{n+m-1}: z^{\bz/l}({\b A}^{n+m-1})^{\mu_l}\to z^{\bz/l}({\b A}^{n+m-1})^{h\mu_l}$$ 
induces isomorphisms on the factors $K(\bz/l,2i,i)$ for $n\leq i\leq n+m-1$. 

Consider the following diagram
$$\xymatrix{
z^{\bz/l}({\b A}^{n+m-1})^{\mu_l}\ar[r]^{\jm} \ar[d]^{\gamma_{n+m-1}} & z^{\bz/l}({\b A}^{n+m})^{\mu_l}\ar[d]^{\gamma_{n+m}}\\
z^{\bz/l}({\b A}^{n+m-1})^{h\mu_l}\ar[r]^{\jm_h} & z^{\bz/l}({\b A}^{n+m})^{h\mu_l},
}$$
where $\jm$ and $\jm_h$ are induced by the closed embedding ${\b A}^{n+m-1}\to {\b A}^{n+m}$ with the last coordinate zero. The diagram is clearly commutative by our construction of $\gamma$. 

We have proved that $\gamma_{n+m}$ induces the identity on the factor $K(\bz/l,2(n+m),n+m)$ by Lemma \ref{ltop}, so by the induction hypothesis we only have to show that $\jm$ and $\jm_h$ induce isomorphisms on the factors $K(\bz/l,2i,i)$ for $n\leq i\leq n+m-1$. 

This is clearly true for $\jm$ by how we compute the type of $z^{\bz/l}({\b A}^{n+m})^{\mu_l}$ (see the proof Theorem \ref{pz}, especially \er{lowerones} and \er{gl}). 

Now consider $\jm_h: z^{\bz/l}({\b A}^{n+m-1})^{h\mu_l}\to z^{\bz/l}({\b A}^{n+m})^{h\mu_l}$. For a scheme $X$, we have the following commutative diagram
\bean
\xymatrix{
Hom_{H^{\ba^1}_\bullet}(X_+,z^{\bz/l}({\b A}^{n+m-1})^{h\mu_l})\ar[r]^\cong\ar[d]^{\jm_h\circ} & H^{2(n+m-1),n+m-1}(X\times B\mu_l,\bz/l)\ar[d]^{\jm_H}\\
Hom_{H^{\ba^1}_\bullet}(X_+,z^{\bz/l}({\b A}^{n+m})^{h\mu_l})\ar[r]^\cong & H^{2(n+m),n+m}(X\times B\mu_l,\bz/l),
}\l{indu}
\eean
where the two isomorphisms are by \er{bcnm}, the left vertical map $\jm_h\circ$ is composition by $\jm_h$. The right vertical map $\jm_H$ is defined as the follow composition:
\bean
&\jm_H:&H^{2(n+m-1),n+m-1}(X\times B\mu_l,\bz/l)\nm\\
&\os {\cdot t_{L_m}}\cong& \widetilde H^{2(n+m),n+m}(X_+\wedge Th(L_m\to B\mu_l),\bz/l)\nm\\
&\os {i^*} \to& H^{2(n+m),n+m}(X\times B\mu_l,\bz/l),\nm
\eean
where $Th(L_m\to B\mu_l)$ is the Thom space of the line bundle $L_m$ on $B\mu_l$ associated to the last ${\b A}^1$ component of ${\b A}^m$, the isomorphism is the Thom isomorphism via the multiplication by the Thom class $t_{L_m}$, and the last arrow is induced by inclusion $i:B\mu_l\to Th(L_m\to B\mu_l)$ as the zero section. The commutativity of the diagram \er{indu} is seen from the proof of \er{bcnm} (see \cite[Section 5]v), where one adds a vector bundle to make the direct sum trivial and applies the Thom isomorphism theorem to pull back the dimension of the cohomology class. 

Therefore 
$$\jm_H=\cdot e(L_m): H^{2(n+m-1),n+m-1}(X\times B\mu_l,\bz/l)\to H^{2(n+m),n+m}(X\times B\mu_l,\bz/l)$$
is the multiplication by $e(L_m)$, the Euler class of the line bundle $L_m$.
Since the $\mu_l$ action on ${\b A}^m$ is a direct sum of nontrivial irreducible representations, we have $e(L_m)=a_m v$, where $1\leq a_m\leq l-1$ is an integer invertible in $\bz/l$ and $v$ is the Euler class of the line bundle on $B\mu_l$ associated to the standard representation $\rho$, which by definition is the $v$ in \er{hnm}. 


Therefore in view of \er{hnm}, $\jm_H$ has the form 
\bea
&\jm_H:&H^{2(n+m-1),n+m-1}(X\times B\mu_l, \bz/l)\cong \bigoplus_{\epsilon=0,1,\ 0\leq \delta\leq n+m-1} H^{*,*}(X,\bz/l)u^\epsilon v^\delta\\
&\os {\cdot a_m v} \lra& \bigoplus_{\epsilon=0,1,\ 0\leq \delta\leq n+m-1} a_m H^{*,*}(X,\bz/l)u^\epsilon v^{\delta+1}\to H^{2(n+m),n+m}(X\times B\mu_l,\bz/l),
\eea
which on the components of $H^{p,*}(X,\bz/l)$ for $p\leq 2(n+m-1)$ are multiplications by $a_m$ and thus isomorphisms. This proves our theorem.
\end{pfmn}

Finally let's add some remarks about our operations.

\begin{rmk}  Our construction doesn't produce any cohomology operations lowering algebraic weights. Meanwhile Voevodsky's approach {\it a priori} produces such operations (see \er{hcl}). It takes some effort for Voevodsky (see \cite[Section 3]v) to prove that all such operations lowering algebraic weights are zero. 
\end{rmk}

\begin{rmk}  Voevodsky's operations $H^{2n,n}(X,\bz/l)\to H^{2n+2i,n+i}(X,\bz/l)$ are nonzero only if $i=j(l-1)$ is a multiple of $l-1$. The proof uses the full symmetry of the construction under the whole symmetric group $S_l$. By our comparison theorem \ref{tcomp}, we see that this also holds for our operations, i.e. $D^i=0$ and thus $E^i=0$ by Theorem \ref{pbock} unless $i=j(l-1)$. It would be satisfying to give a direct proof of this fact for our operations. But this eludes my effort so far mainly because in our construction we are dealing with finite dimensional lens spaces, and the homotopy classes of maps between such lens spaces are more complicated than those for infinite dimensional lens spaces, which were employed in Voevodsky's approach. 
\end{rmk}

\begin{rmk}  Using a trick of considering only bistable operations, Voevodsky can extend his operations to motivic cohomology $H^{*,*}(X,\bz/l)$ of all dimensions instead of just $H^{2n,n}(X,\bz/l)$, and furthermore by a formality (see \cite[Corollary 2.10]v) all bistable operations are group homomorphisms. To follow his approach, we need to verify that for us 
\bean \l{stability}
{{\cal P}}_{n+1}(x\wedge \sigma_T)={{\cal P}}_n(x)\wedge \sigma_T,
\eean
where $x\in H^{2n,n}(X,\bz/l)$ is a cohomology class of dimension $(2n,n)$, $\sigma_T\in \widetilde H^{2,1}(T,{\Bbb Z}/l)$ is the tautological motivic cohomology class of $T={\Bbb A}^1/({\Bbb A}^1-\{0\})$, i.e. $\sigma_T$ is represented by the diagonal cycle, and $x\wedge \sigma_T\in \widetilde{H}^{2(n+1),n+1}(X_+\wedge T,\bz/l)$ is the product. Here we write $\cp_n$ to note $\cp$ on dimension $(2n,n)$. This "should" be obtained by the Cartan formula, which "should" say that ${{\cal P}}_{n+1}(x\wedge \sigma_T)={{\cal P}}_{n}(x)\wedge {{\cal P}}_{1}(\sigma_T)$, together with the fact that 
\bean
{{\cal P}}_1(\sigma_T)=\sigma_T.
\l{sigma}
\eean
(\ref{sigma}) holds since the fiber product of the diagonal cycle (being a one-to-one cycle) is supported on the diagonal $\ba^1$ of $\ba^{l}$ and is restricted to zero off the diagonal to $\ba^{l}-\ba^1$, and so ${\cal P}_1(\sigma_T)=D^0(\sigma_T)=\sigma_T$ by Theorem \ref{p0}. But the full powered Cartan formula surely is out of our reach, since it involves some coefficients coming from the base field when $l=2$ (see \cite[Proposition 9.7]v). I have a primitive version of the Cartan formula, but it is not clear to me if it is enough for the stability condition \er{stability}.
\end{rmk}

\begin{rmk}  I want to say, as my last words, that although our construction can be said to be conceptually simpler, we lose much control of the actual classes, especially about the ring structure. In particular, our computational ability is very limited using this approach, and for instance, we don't touch the Adem relations at all (see \cite[Section 10]v). 
\end{rmk}

Department of Mathematics, Texas A\&M University, College Station, TX 77843-3368.

{\it E-mail address:} {\tt nie@math.tamu.edu}


\begin{thebibliography}{99}


\bibitem{dv} Deligne, Pierre, {\it Voevodsky's lectures on motivic cohomology.}  www.math.ias.edu/$\tilde\ $vladimir/seminar.html, 2000/2001.






\bibitem{fv} Friedlander, Eric M.; Voevodsky, Vladimir, 
{\it Bivariant cycle cohomology.}
Cycles, transfers, and motivic homology theories, 138--187,
Ann. of Math. Stud., 143,
Princeton Univ. Press, Princeton, NJ, 2000. 



\bibitem{k0}  Karoubi, Max, {\it Formes diff\'erentielles non commutatives et op\'erations de Steenrod.} (French) [Noncommutative differential forms and Steenrod operations]  Topology  34  (1995),  no. 3, 699--715.

\bibitem{k} Karoubi, Max, 
{\it Produit cyclique d'espaces et op\'erations de Steenrod.} (French) [Cyclic product of spaces and Steenrod operations]
The Arnold-Gelfand mathematical seminars, 281--323,
Birkh\"auser Boston, Boston, MA, 1997. 



\bibitem{mv} Morel, Fabien; Voevodsky, Vladimir, {\it  $\ba^1$-homotopy theory of schemes.}  Inst. Hautes \'Etudes Sci. Publ. Math.  No. 90 (1999), 45--143 (2001).

\bibitem{flat} Raynaud, Michel; Gruson, Laurent, 
{\it Crit\'eres de platitude et de projectivit\'e. Techniques de "platification" d'un module.} (French)
Invent. Math. 13 (1971), 1--89.

\bibitem{s} Steenrod, N. E., {\it Cohomology operations.} Lectures by N. E. Steenrod written and revised by D. B. A. Epstein. Annals of Mathematics Studies, No. 50 Princeton University Press, Princeton, N.J. 1962 vii+139 pp. 


\bibitem{sv2}  Suslin, Andrei; Voevodsky, Vladimir, {\it Relative cycles and Chow sheaves.}  Cycles, transfers, and motivic homology theories,  10--86, Ann. of Math. Stud., 143, Princeton Univ. Press, Princeton, NJ, 2000. 

\bibitem{icm98} Voevodsky, Vladimir, {\it  $\ba^1$-homotopy theory.} Proceedings of the International Congress of Mathematicians, Vol. I (Berlin, 1998).  Doc. Math.  1998,  Extra Vol. I, 579--604 (electronic). 

\bibitem{dm} Voevodsky, Vladimir, {\it Triangulated categories of motives over a field.}  Cycles, transfers, and motivic homology theories,  188--238, Ann. of Math. Stud., 143, Princeton Univ. Press, Princeton, NJ, 2000.

\bibitem{allchar}  Voevodsky, Vladimir, {\it Motivic cohomology groups are isomorphic to higher Chow groups in any characteristic.}  Int. Math. Res. Not.  2002,  no. 7, 351--355.


\bibitem{v}  Voevodsky, Vladimir, {\it Reduced power operations in motivic cohomology.}  Publ. Math. Inst. Hautes \'Etudes Sci.  No. 98 (2003), 1--57. 


\end{thebibliography}
\end{document}